\documentclass[12pt,reqno]{amsart}
\usepackage[letterpaper, margin=1.2in]{geometry}
\usepackage{amstext,amsthm,amsmath,amsfonts,amssymb,amscd,latexsym,txfonts,stmaryrd,enumerate}
\usepackage[colorlinks=true, linkcolor=blue, anchorcolor=black, citecolor=blue, filecolor=blue, menucolor= blue, urlcolor=black]{hyperref}
\usepackage{graphicx}
\usepackage{epstopdf}
\usepackage{hyperref}
\usepackage{setspace}
\usepackage[all]{xy}
\usepackage{calrsfs}     

\usepackage{tikz}
\usepackage{tikz}
\usepackage{pgflibraryarrows}
\usepackage{pgflibrarysnakes}

\def\-{\raisebox{.75pt}{-}}

\newtheorem{theorem}{Theorem}[section]

\newtheorem{cor}[theorem]{Corollary}

\newtheorem{alg}[theorem]{Algorithm}
\newtheorem{proposition}[theorem]{Proposition}
\newtheorem{lemma}[theorem]{Lemma}

\theoremstyle{definition}
\newtheorem{definition}[theorem]{Definition}

\newtheorem{remark}[theorem]{Remark}

\newtheorem{example}[theorem]{Example}

\newcommand{\ifi}{\ensuremath{\ \Leftrightarrow}\ }

\newcommand{\reg}{\ensuremath{\mathrm{reg}}}
\newcommand{\DReg}{\ensuremath{\mathrm{DReg}}}
\newcommand{\Ann}{\ensuremath{\mathrm{Ann}}}
\newcommand{\Ass}{\ensuremath{\mathrm{Ass}}}
\newcommand{\Min}{\ensuremath{\mathrm{Min}}}
\newcommand{\coker}{\ensuremath{\mathrm{coker}}}
\newcommand{\Det}{\ensuremath{\mathrm{Det}}}
\newcommand{\h}{\ensuremath{\mathrm{ht}}}
\newcommand{\img}{\ensuremath{\mathrm{im}}}

\newcommand{\het}{\ensuremath{\mathrm{ht}}}

\newcommand{\fp}{\ensuremath{\mathfrak{p}}}
\newcommand{\fq}{\ensuremath{\mathfrak{q}}}
\newcommand{\pp}{\ensuremath{\mathbb{P}}}

\newcommand{\Res}{\ensuremath{\mathrm{Res}}}
\DeclareMathOperator{\Spec}{Spec}
\DeclareMathOperator{\St}{St}
\newcommand{\Sym}{\ensuremath{\mathrm{Sym}}}

\newcommand{\Q}{\ensuremath{\mathbb{P}_k^1\times\mathbb{P}_k^1}}
\renewcommand{\P}{\ensuremath{\mathbb{P}_k}}
\newcommand{\N}{\ensuremath{\mathbb{N}}}
\newcommand{\Z}{\ensuremath{\mathbb{Z}}}

\newcommand{\cB}{\ensuremath{\mathcal{B}}}
\newcommand{\cM}{\ensuremath{\mathcal{M}}}
\newcommand{\cR}{\ensuremath{\mathcal{R}}}
\newcommand{\be}{\ensuremath{\mathbf{c}}}
\newcommand{\F}{\ensuremath{\mathbb{F}}}
\newcommand{\B}{\ensuremath{\mathbb{B}}}

\title{ Implicitization of tensor product surfaces via virtual projective resolutions}

\author{Eliana Duarte, Alexandra Seceleanu}
\date{\small \today}                                           
\address{Eliana Duarte \\ Max-Planck-Institute for Mathematics in the Sciences, Leipzig and Otto-von-Guericke Universit\"at 
, Magdeburg \\ eliana.duarte@ovgu.de}
\address{Alexandra Seceleanu \\ Mathematics Department\\University of Nebraska--Lincoln\\ Lincoln, NE 68588\\ aseceleanu@unl.edu}
\keywords{implicitization, resultant, residual resultant, virtual projective resolution, multigraded regularity, Eagon-Northcott 
complex.\\
\textit{2010 Mathematics Subject Classification}: Primary: 13P15;  Secondary: 13D02, 14Q10.}

\begin{document}
\maketitle

\begin{abstract}
We derive the implicit equations for certain parametric surfaces in three-dimensional projective space termed tensor product surfaces. Our method computes the implicit equation for such a surface based on the knowledge of the syzygies of the base point locus of the parametrization by means of constructing an explicit virtual projective resolution. 
\end{abstract}

\section{\large Introduction}

The \emph{residual resultant} of a system of polynomial equations is a polynomial on the coefficients of the system 
that vanishes if and only if the system has a solution outside the zero set of another prescribed system of polynomial equations. Residual resultants for 
projective space were introduced in \cite{residualRes} and further developed in \cite{residualBuse} for the case of $\pp^2$. In this article we consider residual resultants over $\Q$. 

For projective space, the computation of the residual resultant relies on producing a free resolution of
an ideal having the same vanishing locus as the residual (colon) ideal of the two systems of polynomial equations. In this article we formulate a similar  approach to compute a residual resultant over $\Q$ where we replace the free resolution of the residual ideal with a \emph{virtual
resolution}. This allows the derivation of the residual resultant from smaller, more manageable complexes 
than the more standard free resolutions. Besides being shorter than their free resolution counterparts, virtual resolutions also exhibit a closer relationship with Castelnuovo--Mumford regularity than minimal free resolutions.
We exploit this relationship and present Algorithm \ref{alg:residualres} to compute residual resultants
over $\Q$. 

Our motivation to study residual resultants over $\Q$ comes from implicitization in geometric modeling. 
In this context,  a {\em tensor product surface}  is the closure of the image $\Lambda$ of a rational map $\lambda: \Q \dashrightarrow \P^3$  defined by four bihomogeneous polynomials $p_0, p_1, p_2, p_3\in k\left[\Q\right]=k[s,t,u,v]$ as
$$\lambda([s:t],[u:v])=[p_0(s,t,u,v):p_1(s,t,u,v):p_2(s,t,u,v):p_3(s,t,u,v)].$$ 
The {\em base points} of $\lambda$ are the common zeros of the polynomials $p_0, p_1, p_2, p_3$.
The {\em implicitization problem} for tensor product surfaces consists on finding the equation whose vanishing 
defines the surface $\Lambda$ in $\pp^3$. This problem has its origins in the seminal papers \cite{SederbergChen,CoxGoldmanZhang} and has been considered further in \cite{GSK,D'Andrea,Botbol}.

Three methods can be used to solve the implicitization problem for tensor
product surfaces: Gr\"obner bases, resultants, and Rees algebras. 
Gr\"obner basis methods are least satisfactory since they  tend to be computationally intensive. Thus, it is primarily the latter two techniques which are used. Since classical resultants fail in the presence of base points,
 following the work of Bus\'e \cite{residualBuse},  we propose the use of residual resultants over $\Q$ to solve the
implicitization problem for tensor product surfaces in this case. We present this approach in Algorithm~\ref{alg:implicitization}.

 The structure of this paper is as follows: in section \ref{s:resultant} we give the necessary background on residual resultants, with special attention to the case of biprojective space. In section \ref{s:compresultant} we derive effective methods to compute the residual resultant based on a virtual projective resolution for certain ideals of minors. In section \ref{s:implicitization} we show how this theory can be applied to the implicitization problem for tensor product surfaces. Finally, section \ref{s:examples} contains many worked out examples that illustrate our results. 

Throughout the paper $\N$ denotes the set of nonnegative integers. 
 

\section{\large A residual resultant for \Q}\label{s:resultant}

In this section we give an overview of the theory and construction for a residual resultant over a biprojective space. We follow closely the exposition in \cite{residualRes} and \cite{residualBuse} adapting the statements for the case of the variety $Q= \Q$. 

Algebraically, classical resultant computations can be phrased as follows: 
given  commutative rings $A=k[x_0,\ldots, x_m]=k[\pp^m_{k}]$ and $C=k[C_{ij} : 0\leq i \leq n, 1\leq j \leq \dim_k(A_{d_i})]$, where the latter is viewed as a ring of indeterminate coefficients,  form the polynomial ring $T=C[x_0,\ldots, x_m]=C\otimes_kA$ and define a set of homogeneous polynomials $F_0,\ldots, F_m\in T$
\[
F_i(C_{ij}, x_0,\ldots, x_n)=\sum_{m_j\in A_{d_i}} C_{ij}m_j.
\]
One is interested in finding a generator for the principal ideal $I=(F_0,\ldots, F_m)\cap C$, which is called the resultant of $F_0,\ldots F_m$. The resultant is a unique (up to scaling by constants) irreducible polynomial in $C$ \cite[Chapter 12]{GKZ}.  For a point $c=(c_{ij}) \in \mathbb{P}_k(C_1)$ define the evaluation map at $c$ to be the $A$-module homomorphism $e_c:T\to A,  e_c(C_{ij})=c_{ij}$ induced by the analogous $k$-linear map $ C \to k(c)$. 
The zero locus of the ideal $I$
$$V(I)=\{c \in \mathbb{P}_k(C_1) : V\left(e_c(F_0),\ldots, e_c(F_m)\right)\neq \emptyset \}$$
consists of the coefficients $c=(c_{ij})$ for which the equations $e_c(F_0),\ldots, e_c(F_m)$ have common solutions in $\P^m$. 

We proceed to describe a modified version of this classical resultant termed the residual resultant. If $A$ is the coordinate ring of a variety $Q$ and $C,T$ are as above, consider 
 two sets of homogeneous polynomials $g_0,\ldots, g_n\in A$ and 
 $F_0,\ldots  F_m\in (g_0,\ldots, g_n)T$. 
The residual resultant is a generator for the principal ideal $I=(F:G)\cap C$, where $F=(F_0,\ldots, F_m)$ and $G=(g_0,\ldots, g_n)$. The zero locus of this ideal
\[
V(I)=\{c \in \mathbb{P}_k(C_1): V\left(e_c(F_0),\ldots, e_c(F_m) \right) \setminus V\left(g_0,\ldots, g_n \right) \neq \emptyset \}
\]
consists of the coefficients $c_{ij}$ for which the equations $e_c(F_0),\ldots, e_c(F_m)$ have common solutions outside the common zero locus of $g_0,\ldots g_n$ in $Q$.

We now rephrase the problem in the language of algebraic geometry. The classical resultant is interpreted in this language in \cite{Jou1, Jou2} and \cite[Propositions 3.1 and 3.3]{GKZ}. Following the exposition in \cite{residualRes}, let $Q$ be a an irreducible projective variety of dimension $\dim(Q)=m$ over the algebraically closed field $k$. 
Consider $m+1$ invertible sheaves $\mathcal{L}_0,\ldots, \mathcal{L}_m$ on $Q$ and let $V_i=H^0(Q,\mathcal{L}_i)$ be the vector space spanned by  the global sections
of the sheaf $\mathcal{L}_i$.  Poposition~\ref{prop:1} sets up the residual resultant as
a polynomial that captures the condition for a set of global sections $f_0,\ldots,f_m$ ($f_i \in V_i$) to vanish on the variety $Q$. This resultant is a polynomial in the coefficients of each $f_i$ with respect to the basis of the vector space $V_i$. 

More precisely, given a set of polynomials  $F_i=\sum_{b_j\in B_i} C_{ij}b_j\in T$ expressed in terms of fixed bases $B_i$ for each vector space $V_i$, their resultant is a polynomial $\Res_{V_0,\ldots,V_m}\in C$. For any $c\in \P(C_1)$, if $f_i=e_c(F_i)$, then  $\Res_{V_0,\ldots,V_m}(f_0,\ldots, f_m)$ denotes the polynomial 
\[
\Res_{V_0,\ldots,V_m}(f_0,\ldots, f_m) = e_c\left(\Res_{V_0,\ldots,V_m}\right).
\]
From this point onward, we use the notation $F_0,\ldots,F_m$ for elements of $T$ and  $f_0,\ldots, f_m$ for  specializations  $f_i=e_c(F_i)$ at some $c\in \P(C_1)$.

\begin{proposition}[{\cite[Proposition 1]{residualRes}}]\label{prop:1}
Suppose that each $V_i$ generates the sheaf $\mathcal{L}_i$ on $Q$ and that $V_i$ is very ample on a nonempty open subset $U$ of $Q$. Then there exists an irreducible polynomial on $\prod\limits_{i=0}^mV_i$, denoted by $\Res_{V_0,\ldots,V_m}$ and called the $(V_0,\ldots,V_m)$-resultant, which satisfies 
\[
\Res_{V_0,\ldots,V_m}(f_0,\ldots, f_m)= 0 \iff \exists \ x \in X : f_0(x) = \dots = f_m(x) = 0.
\]
Moreover, $\Res_{V_0,\ldots,V_s}$ is homogeneous in the coefficients of each $f_i$, and of
degree $\int_Q \prod_{j\neq i} c_1(\mathcal{L}_j)$.
\end{proposition}

We will follow the aforementioned result to  define a {\em residual resultant} for $Q=\Q$. This follows readily using the methods of \cite{residualRes},  but it is important for our purposes to establish the notation in terms of sheaves on $\Q$ instead of $\P^n$. For this reason we include a discussion of the setup below. 

From this point on let $R=k[s,t,u,v]$ denote the bigraded coordinate ring of $\Q$ over an algebraically closed
field $k$, with $\deg(s)=\deg(t)=(1,0)$ and $\deg(u)=\deg(v)=(0,1)$.  
Let $R_{(a,b)}$ denote the set of elements in $R$ of bidegree $(a,b)$.
Recall that the smallest geometrically irrelevant ideal of $\Q$ is $B=(s,t)\cap(u,v)$. This yields a family of geometrically irrelevant ideals for $\Q$, i.e. $\cB=\{\fp \in \Spec(R): B\subseteq \fp\}$. 

\begin{definition}
The $B$-saturation of an ideal $I\subset R$ is the ideal $I^{\rm sat}=\bigcup_{i=0}^\infty I:B^i$, where
$I:B^i=\{f\in R : fB^i\in I\}$. The  geometric importance of the $B$-saturation stems from the fact  that for bihomogeneous ideals $I\subseteq R$, the following varieties agree $V(I)=V(I^{\rm sat})$. Analogously one defines the $B$-saturation of an $R$-module $M$ to be $M^{\rm sat}=H^0\left(\widetilde M, \Q\right)$.
\end{definition}

Let $Q=\Q=\mathrm{Proj}(R)$ and consider a bihomogeneous ideal $G=(g_1,\ldots,g_n)\subseteq R$ where $\deg g_j= (k_j,l_j)$. Let $\mathcal{G}$ be the coherent sheaf of ideals associated to $G$. Consider pairs of nonnegative integers $(a_i,b_i),\, 0\leq i\leq 2,$ such that $(a_i,b_i)\geq (k_j,l_j)$ entrywise  for all $i,j$, which yield the sheaves $\mathcal{G}(a_i,b_i)=\mathcal{G}\otimes_{\mathcal{O}_Q}\mathcal{O}_{Q}(a_i,b_i)$ for $0\leq i\leq 2$. The vector space $V_i=H^0(Q,\mathcal{G}(a_i,b_i))$ is the set of polynomials of degree $(a_i,b_i)$ which belong to the saturation of the ideal $G$. We denote by $\pi:\widetilde{Q}\to Q$ the blow-up of $Q$ along the sheaf of ideals $\mathcal{G}$. The inverse image of the sheaf $\widetilde{\mathcal{G}}=\pi^{-1}\mathcal{G}\cdot \mathcal{O}_Q$ is an invertible sheaf on $\widetilde{Q}$. The sheaf $\widetilde{\mathcal{G}}\otimes \pi^\ast(\mathcal{O}_{Q}(a_i,b_i))$ is denoted by $\widetilde{\mathcal{G}}(a_i,b_i)$.

Proposition \ref{prop:tpres} establishes the existence of a residual resultant polynomial, which cuts out the locus of those polynomials $f_0,f_1, f_2\in V_0\times V_1\times V_2$ for which the common vanishing  of $f_0,f_1,f_2$ contains a point not in $V(G)$. It also gives an algebraic criterion for this geometric condition, namely that the saturations of the two ideals $G=(g_1, \dots, g_n)$ and $F=(f_0,f_1,f_2)$ with respect to $B$ are distinct. In order to establish this fact we need the following definition.

\begin{definition}
An ideal $I\subseteq R$ is said to be {\em locally a complete intersection} if $I_\fp$ can be generated by a regular sequence for every prime ideal $\fp\in \Spec(R)\setminus \cB$.
\end{definition}

\begin{proposition}[{\cite[Proposition 3]{residualRes}}] \label{prop:tpres}
Let $G=(g_1,\ldots,g_n)\subseteq R$ be a codimension two locally complete intersection ideal, with $\deg(g_j)=(k_j,l_j)$. Choose bihomogeneous polynomials $f_i\in V_i=\mathcal{G}(a_i,b_i)$ for $i=0,1,2$ such that $F=(f_0,f_1,f_2)$   and the following condition holds
 \[
 (a_i,b_i)\geq (k_{j_1}+1,l_{j_1}) \text{ for some } j_1 \text{ and } (a_i,b_i)\geq (k_{j_2},l_{j_2}+1) \text{ for some } j_2.
 \]
Then there exists a 
polynomial in $C=\prod_{i=0}^{2}k\left[V_i\right]$, denoted $\mathrm{Res}_{\mathcal{G},\{(a_{i},
b_{i})\}_{i=0}^2}$ which satisfies
\begin{eqnarray}
\mathrm{Res}_{\mathcal{G},\{(a_{i},b_{i})\}_{i=0}^2}(f_0,f_1,f_2)=0 &\ifi& \exists \ x\in \widetilde{Q}:\pi^{\ast}(f_{0})(x)=\pi^{\ast}(f_{1})(x)=\pi^{\ast}(f_{2})(x)=0 \\
&\ifi& \exists \ y\in \Q \text{ such that } y\in V(F)\setminus V(G) \\
&\ifi&F^{\rm sat}\neq G^{\rm sat}.
\end{eqnarray}
\end{proposition}

\begin{proof} Let $i\in\{0,1,2\}$  and consider the vector space of global sections $V_i=H^0(Q,\mathcal{G}(a_i,b_i))$.
The sections $s\in V_i$ generate the invertible sheaf $\mathcal{G}(a_i,b_i)$ on an open subset of $Q$, namely $Q\setminus Z$.
Following \cite[Ch.II.7.17.3]{Hartshorne} we blow-up $\Q$ at the subscheme defined by $\mathcal{G}$.
Then $\widetilde{\mathcal{G}}(a_i,b_i)$ is  globally generated by the pullbacks $\pi^*(s)$ for 
 $s\in H^0(Q,\mathcal{G}(a_i,b_i))$.  Thus for all $i\in\{0,1,2\}$, if we
let $\widetilde{V}_i$ be the vector subspace generated by the pullbacks $\pi^*(s),\, s\in V_i$ then $\widetilde{V}_i$ generates $\widetilde{\mathcal{G}}(a_i,b_i)$ on $\widetilde{Q}$. 

Next we show that each $\widetilde{\mathcal{G}}(a_i,b_i)$ is very ample on an open subset $U$ of $\widetilde{Q}$. Suppose $(a_i,b_i)$ satisfies the inequality conditions in
the statement of the proposition. Let $S_k$ be the subvariety of $\widetilde{Q}$ defined by the vanishing of $\pi^{*}(g_{j_k})$ and let $U_k=\widetilde{Q}\setminus S_k$ for $k=1,2$. Set $U= U_1 \cap U_2$. We show that the map $\Gamma_i:U\to\pp(\widetilde{V}_i),\, 
x\mapsto \{\pi^*(f) \mid f\in V_i, \pi^*(f)(x)=0\}$ is an embedding. 
Since a point in \Q\ is a pair $(p_1, p_2)$ where $p_i$ are 
points in the $i$-th factor, there is a form $L_1$ of bidegree $(1,0)$ or $L_2$ of bidegree $(0,1)$ that vanishes at the given 
point but not at another point $(q_1,q_2)\in \Q$ according to whether $p_1\neq q_1$ or $p_2\neq q_2$. We say that such a form 
separates $(p_1, p_2),(q_1, q_2)$.  In the former case there is a global section in $\widetilde{V_i}$ which is a multiple of 
$L_1g_{j_1}$ and which separates $\pi^*(p_1,p_2)$ and $\pi^*(q_1,q_2)$ in $U$. Analogously, if $(p_1, p_2),(q_1, q_2)$ are separated by 
a form of bidegree $(0,1)$, there is a global section in $\widetilde{V_i}$ which is a multiple of $L_2g_{j_2}$ and which 
separates $\pi^*(p_1,p_2)$ and $\pi^*(q_1,q_2)$ in $U$. A proof that the differential condition for very ampleness holds follows in a similar fashion to \cite[Proposition 3]{residualRes} by the use of the appropriate separating form in each case. We conclude that each $\widetilde{\mathcal{G}}
(a_i,b_i)$ is very ample on the non-empty open subset $U$.

The first equivalence of the conclusion follows by applying Proposition~\ref{prop:1} to the invertible sheaves $\widetilde{\mathcal{G}}(a_i,b_i)$ on $\widetilde{Q}$. For $(2)\Rightarrow (1)$ notice that if $y\not\in V(G)$ and $f_0(y)=f_1(y)=f_2(y)=0$ then, for the  unique $x\in\widetilde{Q}$ such that $\pi(x)=y$, we have $\pi^{\ast}(f_{0})(x)=\pi^{\ast}(f_{1})(x) = \pi^{\ast}(f_{2})(x)=0$.  The equivalence $(2)\Leftrightarrow (3)$ follows from the identities $V(F)=V(F^{\rm sat})$ and $V(G)=V(G^{\rm sat})$. 
It remains to show that $(3)\Rightarrow (1)$, equivalently, if $\pi^*(f_0),\,\pi^*(f_1),\,\pi^*(f_2)$ do not vanish simultaneously on $\widetilde{Q}$ then $F^{\rm sat}= G^{\rm sat}$. Since $G$ is locally a complete intersection, the sheaf $\mathcal{G}/\mathcal{G}^2$ is locally free of rank $2$. 
%
 Hence, setting $\mathcal{F}$ to be the ideal sheaf corresponding to $F$, one sees that the inclusion $\mathcal{F}\hookrightarrow \mathcal{G}$ is a surjection locally at $\fp\in X$. Thus $\mathcal{F}=\mathcal{G}$ and hence $F^{\rm sat}=G^{\rm sat}$ holds true.
\end{proof}

\begin{remark}
By the assumption on the codimension of $G$, the ideal sheaf $\mathcal{G}$ in Proposition \ref{prop:tpres} defines a zero dimensional scheme. Proposition \ref{prop:tpres} applies when $G$ defines a reduced set of points in \Q, since such an ideal is locally a complete intersection  by \cite[Lemma 4.1]{CFGLMNSSV}. However, not all ideals $G$ that fit the hypotheses of Proposition \ref{prop:tpres} define reduced sets of points in \Q. For example $G=\langle s^2t^2, u^2v^2 \rangle$ is a (global) complete intersection, hence this ideal is also locally a complete intersection, which is not reduced.
\end{remark}

Suppose that the ideal sheaf $\mathcal{G}$ defines a zero dimensional scheme $Z$ composed of $p$ points $P_1,\ldots,P_p$. We denote by $e_i$ the the multiplicity
of the point $P_i$ in $Z$. We have
\[e_i=\dim_{k}(\mathcal{O}_{Z,P_i}), \mbox{ where } \mathcal{O}_{Z}=\mathcal{O}_{\Q}/\mathcal{G},\]
and hence $\sum_{i=1}^{p}e_{i}=\dim_{k}H^{0}(Z,\mathcal{O}_{Z})$. 

\begin{remark} 
\label{rem:incidence}
One important aspect to recall from the proof of Proposition~\ref{prop:1} \cite[Proposition 1]{residualRes} is that the incidence variety
defined by 
\[
\widetilde{W}= \left\{(x,f_0,\ldots,f_m)\in \widetilde{Q}\times \prod_{i=0}^{m}\pp(V_i) :  f_0(x)=\ldots=f_m(x)=0\right\} \subseteq \widetilde{Q}\times \prod_{i=0}^{m}\pp(V_i)
\]
 has codimension $m+1$.
In the context of Proposition~\ref{prop:tpres}, $m=\dim \widetilde{Q}=2$ because $\widetilde{Q}$ is the blowup of $\Q$
at the scheme $Z$ defined by  $\mathcal{G}$.
 Therefore the incidence variety $\widetilde{W}$
in this case is contained in 
$\widetilde{Q}\times \prod_{i=0}^{2}\pp(V_i)$ and it is of codimension $3$. Let $E$ denote the exceptional locus of the blow-up of $\Q$ at $Z$. Then $\widetilde{Q}\setminus E$ is isomorphic to $Q\setminus Z$. The open set 
\[
U=\left\{(x,f_0,f_1,f_2)\in \widetilde{Q}\setminus E\times \prod_{i=0}^{2}\pp(V_i): \pi^*(f_0)(x)=\pi^*(f_1)(x)=\pi^*(f_2)(x)=0\right\}
\]
is dense in $\widetilde{W}$ and  isomorphic to 
\[
W =\left \{(x,f_0,f_1,f_2)\in \left(Q\setminus Z \right)\times \prod_{i=0}^{2}\pp(V_i) :  f_0(x)=f_1(x)=f_2(x)=0\right\}
\] thus  $W$ is of codimension 
three in  $\left(Q\setminus Z\right) \times \prod_{i=0}^{2}\pp(V_i)$.
\end{remark}

In the next proposition we compute the degree of the residual resultant in the
coefficients of each polynomial $f_{i}$. A general formula for this degree is given in 
Proposition~\ref{prop:1} {\cite[Proposition 1]{residualRes}}  and the case for $\pp^{2}$ is 
treated in \cite{residualBuse}. We will now deduce this degree for the residual resultant in $\Q$; the proof follows the same lines as for $\pp^{2}$, except that the computation of the intersection product is now performed on the blow-up of  $\Q$ at $Z$.

\begin{proposition}
\label{prop:degrees}
The polynomial $\mathrm{Res}_{\mathcal{G},\{(a_{i},
b_{i})\}_{i=0}^2}$ is multihomogeneous in the coefficients of each $V_{i}$, of degree $N_{i}$ for $i=0,1,2$ with
\begin{align*}
N_{0} &=  a_{1}b_{2} + b_{1}a_{2} - \sum_{i=1}^{p}e_{i}\,,\, &
N_{1} &=  a_{0}b_{2}+b_{0}a_{2} - \sum_{i=1}^{p}e_{i}\,,\, \text{ and }  &
N_{2} &=  a_{0}b_{1}+b_{0}a_{1} - \sum_{i=1}^{p}e_{i}.
\end{align*}
\end{proposition}
\begin{proof}
 We compute the integer $N_{0}$, the computation of $N_1,N_2$ is carried out in a similar fashion. Fix $i=0$. By Propositions \ref{prop:1}, $N_0$  equals
\[\int_{\widetilde{Q}}c_{1}(\widetilde{\mathcal{G}}_{(a_{1},b_{1})})c_{1}(\widetilde{\mathcal{G}}_{(a_{2},b_{2})})\]
where $c_{1}(\mathcal{F})$ denotes the first Chern class of the sheaf $\mathcal{F}$  over $\widetilde{Q}$ and
$\int_{\widetilde{Q}}$ denotes the degree map on $\widetilde{Q}$. Denote by $H=\pi^{\ast}(h)$ and $L=\pi^{\ast}(l)$ the pullbacks of generic hyperplanes
in $\Q$ that generate the divisor class group $\rm{Cl}(\Q)\cong \mathbb{Z}^{2}$. Each $E_{i}$, $i=1,\ldots,p$ denotes the exceptional divisor of the blow-up $\pi$
above each point $P_{i}$ defined by $\mathcal{G}$, and $E_{i}^{red}$
 the reduced scheme of $E_{i}$. Following \cite{fulton},
$c_{1}(\widetilde{\mathcal{G}}_{(a_{i},b_{i})})=a_{i}H+b_{i}L-\sum_{i=1}^{p}E_{i}.$
Since $E_{i}\cdot E_{j}=0$ if $i\neq j$, $H\cdot E_{i}=L\cdot E_{i}=0$ and $L^{2}=H^{2}=0$, we obtain
\begin{eqnarray*}
\int_{\widetilde{Q}}c_{1}(\widetilde{\mathcal{G}}_{(a_{1},b_{1})})c_{1}(\widetilde{\mathcal{G}}_{(a_{2},b_{2})}) & = & 
                    \int_{\widetilde{Q}}(a_{1}H+b_{1}L-\sum_{i=1}^{p}E_{i})(a_{2}H+b_{2}L-\sum_{i=1}^{p}E_{i})\\
                    &=& \int_{\widetilde{Q}} a_{1}b_{2}H\cdot L +a_{2}b_{1}H\cdot L + \sum_{i=1}^{p}E_{i}^{2}.
\end{eqnarray*}
Now let $f_1$ (resp. $f_2$) be generic global sections of $\mathcal{G}(a_1,b_1)$ (resp. $\mathcal{G}(a_2,b_2)$) and let
$D_{f_1}:=V(f_1)$ (resp. $D_{f_2}:= V(f_2)$) be the divisor corresponding to the vanishing of the section $f_1$ (resp. $f_2$) in $\Q$. We have
\begin{align*}
\pi^{*}D_{f_1} &= \tilde{D}_{f_1}+ \sum_{i=1}^{p} E_i= \tilde{D}_{f_1}+ \sum_{i=1}^{p} m_i E_{i}^{red}, \text{ and }\\
\pi^{*}D_{f_2} &= \tilde{D}_{f_2}+ \sum_{i=1}^{p} E_i= \tilde{D}_{f_2}+ \sum_{i=1}^{p} n_i E_{i}^{red}.
\end{align*}
Where $\tilde{D}_{f_1}$ (resp. $\tilde{D}_{f_2}$) is the strict transform of $D_{f_1}$ (resp. $D_{f_2}$)
and where $m_i$ (resp. $n_i$) is the multiplicity of $f_1$ (resp. $f_2$) at the point $P_i$ \cite[Section 4.3]{fulton}.
Now $ \tilde{D}_{f_1} \cdot \tilde{D}_{f_2}=0$ and since $\mathcal{G}$ is a local complete intersection, for each
poin $P_i\in Z$ we have $m_i n_i=e_i$ \cite[Section 12.4]{fulton}. We deduce that 
\begin{align*}
\sum_{i=1}^{p}E_i^{2} = \sum_{i=1}^{p}m_i n_i\, {E_i^{red}}^2= \sum_{i=1}^{p}e_ i \,{E_i^{red}}^2.
\end{align*}
By the projection formulae, we know that $\int_{\widetilde{Q}}H\cdot L=1$ and $\int_{\widetilde{Q}}{E_i^{red}}^2=-1$. 
Therefore
 \[N_{0}  =  a_{1}b_{2} + b_{1}a_{2} - \sum_{i=1}^{p}e_{i}\]
\end{proof}

We shall give a method for the effective computation of the residual resultant on $\Q$ in section \ref{s:resultantcomputation} after reviewing the notion of virtual complexes on $\Q$, which will prove useful in computing the residual resultants.


\section{\large Virtual resolutions in $\Q$ and multigraded regularity}
\label{s:virtual}

Free resolutions have played an important role in the effective computation of resultants. It is shown in \cite{GKZ} that
the classic projective resultant in $\pp^n$ can be computed via a Koszul complex. In a similar manner, 
\cite{residualBuse} and \cite{residualRes} use the Eagon-Northcott and variants of it to compute residual resultants
with respect to locally complete intersection  ideals  over $\P^2$ and
complete intersection on $\pp^n_{k}$ respectively. The Castelnuovo-Mumford regularity of the ideal resolved by this complex is a crucial ingredient for the computation of the residual resultant and the ability to explicitly exhibit a free resolution has the advantage of giving a straightforward way to calculate the regularity. 
For $\Q$, general recipes for the free resolutions of the analogous ideal are not  available, even under the above mentioned assumptions. We overcome this obstacle by showing that virtual resolutions in $\Q$ have the same good properties that free resolutions have for the computation of  resultants and residual resultants in $\pp^n$ and we give an explicit description for a virtual resolution of certain determinantal ideals.

Two bigraded rings are of central importance for the purpose of residual implicitization on $\Q$. The first is the coordinate ring $R=k[s,t,u,v]$ of $\Q$, equipped with a natural $\Z^2$ grading obtained from viewing $\Z^2$ as the Picard group of $\Q$. For simplicity, we call rings graded by $\Z^2$ bigraded. For a finitely generated $R$-module $M$ and a bidegree $\nu\in \Z^2$, the Hilbert function of $M$ at $\nu$ is $H_M(\nu)=\dim_k M_\nu$.

The second ring of interest is $T=R\otimes_k C$, where $C=k[C_{ij}]$ is a ring of indeterminate coefficients as in section \ref{s:resultant}. Note that $T$ is the coordinate ring of the variety $(\Q)\times \prod_{i=0}^m V_i$ and moreover $\B=BT$ is the irrelevant ideal for this variety.   We equip the ring $T$ with a  $\Z^2$ grading given by $\deg_T(c\otimes r)=\deg_R(r)$ for any $r\in R, c\in C$. Thus $T$ is a finitely generated $C$-algebra with $R_{(0,0)}=C$.  For any bidegree $\nu\in \Z^2$ the bigraded component of $T$ in bidegree $\nu$, $T_\nu= R_\nu \otimes_k C$, is a free $C$-module minimally generated by a basis of $R_\nu$.

\subsection{Virtual resolutions in $\Q$}
Virtual resolutions for $\Q$, also known as $B$-torsion complexes, have been discussed in the literature in \cite{MaclaganSmith} and \cite{CoxDickensteinSchenck} among others. Our interest in these complexes was sparked by \cite{2017virtual}.

An $R$-module $M$ is $B$-torsion if $B^iM=0$ for some $i$.

\begin{definition}
\label{def:virtresR}
A bigraded complex of free $R$-modules $P_i=\bigoplus_j R(-a_{ij},-b_{ij})$ of the form
\[
\mathbf{F}: 0 \longrightarrow P_m \stackrel{\varphi_m}{\longrightarrow}  \cdots  \longrightarrow P_1 \stackrel{\varphi_1}{\longrightarrow} 
 P_0
\]
 is called a {\em virtual resolution} of a module $M$ if $\left(H_0(\mathbf{F})\right)^{\rm sat}\cong M^{\rm sat}$ and all the homology modules $H_i(\mathbf{F})$ with $i>0$ are $B$-torsion.  Note that every  free resolution is automatically a virtual resolution. 
\end{definition}

 Virtual resolutions were introduced in \cite{2017virtual} where it is pointed out how these resolutions  capture the geometry of subvarieties of products of projective spaces in an optimal manner. 
For example, saturated ideals defining finite sets of points in $\pp^2$ have a Hilbert-Burch resolution. This
is not the case for ideals of sets of points in $\Q$, however there is a virtual version of this theorem for points in biprojective space.
\begin{proposition}[{\cite[Corollary 5.2]{2017virtual}}]
 \label{prop:HBpoints}
Every zero-dimensional subscheme $Z$ of $\Q$ has a virtual Hilbert-Burch resolution, i.e., there exists an $(m+1)\times m$ matrix $\varphi$ such that the complex $0 \longrightarrow R^{m+1} \stackrel{\varphi}{\longrightarrow} R^m $ is a resolution for $I_m(\varphi)$ and $V(I_m(\varphi))=Z$. 
\end{proposition}

\begin{cor}
\label{cor:HBpoints}
If $G\subseteq R$ is an ideal defining a not necessarily reduced set of points in $\Q$  there exists an ideal $G'$ such that $G^{\rm sat}=G'^{\rm sat}$, and $G'$ has a Hilbert-Burch resolution. Moreover $G$ is locally a  complete intersection if and only if $G'$ is locally a complete intersection.
\end{cor}
\begin{proof}
The first statement is an algebraic reformulation of Proposition \ref{prop:HBpoints} while the second follows since $G^{\rm sat}=G'^{\rm sat}$ implies that $G_{\fp}=G'_{\fp}$ for $\fp\in \Spec(R)\setminus \cB$. 
\end{proof}

\begin{example} \label{ex:twopoints}
Consider the ideal $I=\langle s,u \rangle \cap \langle t,v\rangle=\langle st,sv,tu,uv\rangle$
of a set of two
 points in $\Q$. A free resolution and a virtual resolution of $I$ with $G=\langle 
sv, tu\rangle$ are shown below. Note how the virtual resolution is much simpler than the free resolution and the ideal $G=\langle sv, tu\rangle$ defines the same variety as $I$.

\[
\xymatrix{ 
0 \ar[r] 
 &
  R \ar[rr]_{\scalebox{0.8}{$\bgroup\begin{pmatrix}v\\
      {-u}\\
      {-t}\\
      s\\
      \end{pmatrix}\egroup$}}
&  & 
 R^4 \ar[rrr]_{\scalebox{0.8}{$
                 \bgroup\begin{pmatrix}{-u}&
     {-v}&
     0&
     0\\
     s&
     0&
     0&
     {-v}\\
     0&
     t&
     {-u}&
     0\\
     0&
     0&
     s&
     t\\
     \end{pmatrix}\egroup$}}
      & & &
R^4 \ar[rrr]_{\scalebox{0.8}{$\bgroup\begin{pmatrix}s t&    t u&  s v&   u v\\   \end{pmatrix}\egroup$}}
 & & &
 R\ar[r]
 &
I \ar[r]
& 0
}
\]

\[
\xymatrix{  0 \ar[r] 
&
R 
\ar[r]_{\scalebox{0.8}{ $\bgroup\begin{pmatrix}{-s v}\\
      t u\\
      \end{pmatrix}\egroup$
}} 
& 
R^2 \ar[r]_{\scalebox{0.8}{$\bgroup\begin{pmatrix} t u&  s v\\   \end{pmatrix}\egroup$}}
& 
G \ar[r]
&
0.
}
\]
\end{example}

This example is an instance of a more general phenomenon.

  \begin{example}
  \label{ex:generalpts}
  If $G$ defines a set $Z$ of $r$ general points in $\Q$, from \cite[Example 5.10]{2017virtual} it follows that $G$ has a virtual resolution 
\[
  \begin{tabular}{cc}
\xymatrix{
 0 \ar[r]& 
 *+\txt{$R(-2,-2p)$}\ar[r]  & 
 *+\txt{$R(-1,-p)^2$}\ar[r] & 
 R 
 }
 & 
 \qquad \text{ if } $r=2p$ \text{ and }\\
\xymatrix{
 0 \ar[r]& 
 *+\txt{$R(-2,-2p-1)$}\ar[r]  & 
 *+\txt{$R(-1,-p)$\\$\oplus$\\$R(-1,-p-1)$}\ar[r] & 
 R
 }
 &
  \qquad \text{ if } $r=2p+1$.
\end{tabular}
\]
In particular, any set of general points in $\Q$ is virtually a complete intersection. Further details on which sets of points in $\Q$ are virtual complete intersections appear in \cite{virtualCI}.
  \end{example}

The notion of virtual resolution can be extended to modules over the ring $T$, where the meaning of the word virtual is understood to be with respect to the irrelevant ideal $BT$. To see why this is a natural extension we start by defining a $T$-module $M$ to be $BT$-torsion if $(BT)^iM=0$ for some $i\geq 0$. The following lemma shows that this notion is equivalent to the notion of $B$-torsion for $R$-modules.

\begin{lemma}
\label{lem:BvsBT}
A $T$-module $M$ is $\B$-torsion if and only if $M$ is $B$-torsion  as an $R$-module.
\end{lemma}
\begin{proof}
Denote by $M_R$ the structure of $M$ as an $R$-module induced by restriction of scalars. The claim follows from the identity $\left(\B^iM\right)_R=B^iM_R$.
\end{proof}
By analogy with Definition~\ref{def:virtresR} we say that a bigraded complex of free $T$-modules $P_i=\bigoplus_j T(-a_{ij},-b_{ij})$ of the form
$
\mathbf{F}: 0 \longrightarrow P_m \longrightarrow  \cdots  \longrightarrow P_1 \longrightarrow P_0
$
 is a {\em virtual resolution} of a $T$-module $M$ if $\left(H_0(\mathbf{F})\right)^{\rm sat}\cong M^{\rm sat}$ and  for $i>0$ the homology modules $H_i(\mathbf{F})$  are $\B$-torsion. In view of Lemma~\ref{lem:BvsBT}, $\mathbf{F}$ is a virtual resolution of the $T$-module $M$ if and only if it is a virtual resolution for the $R$-module $M_R$. 

\subsection{Multigraded regularity: strong and weak forms}
In this paper we make use of a notion of (weak) regularity developed in \cite{MaclaganSmith}. 
Although this applies to modules over a polynomial ring graded by a finitely generated abelian group, we are primarily interested in modules over the rings $R$ and $T$ introduced in the beginning of this section, graded by the group ${\rm Pic}(\Q)=\Z^2=\Z\be_1\oplus \Z\be_2$, where $\be_1,\be_2$ are the standard basis vectors. To explain the notion of bigraded regularity define the sets
\begin{align*}
\St_i &= \begin{cases}
\{(r,s)\in \Z^2 :r+s=-i-1, r<0, s<0\} &\text{ for } i >  0, \\ 
\{(r,s)\in \Z^2 :r+s=-i,  r\geq 0, s \geq 0\} &\text{ for } i \leq  0.
\end{cases}\\
&=\begin{cases}
 \{(-i,-1), (-i+1,-2), \ldots , (-2,-i+1), (-1,-i)\} &\text{ for } i >  0, \\ 
 \{(-i,0), (-i-1,1), ..., (1,-i-1), (0,-i)\} &\text{ for } i \leq  0.
\end{cases}
\end{align*}

\begin{definition}\label{def:wreg}
A module $M$ over a bigraded ring is said to be {\em weakly $\nu$-regular} with respect to the irrelevant ideal  $B$ of that ring if $H_B^i (M)_{\mu} = 0$ for all $i\geq 0$ and $\mu \in \St_i+\nu+\N^2$.  We denote by $\reg(M)$ the set of all elements $\mu\in \Z^2$ such that $M$ is weakly $\mu$-regular and we call this set the {\em regularity region} of $M$.
\end{definition}

As before, the notion of regularity for $T$-modules and $R$-modules are closely related.

\begin{lemma}
\label{lem:BvsBT2}
For a $T$-module $M$ and a bidegree $\mu\in\Z^2$, $M$ is weakly $\mu$-regular with respect to $\B$ if and only if $M$ is weakly $\mu$-regular as an $R$-module with respect to $B$.
\end{lemma}
\begin{proof}
By independence of basis for local cohomology $H^i_B(M)\cong H^i_{\B}(M)$ as $T$-modules, whence $H^i_B(M)_\mu=0$ if and only if $H^i_{\B}(M)_\mu=0$. 
\end{proof}

One of the main applications of (multigraded) regularity consists of controlling the growth of Hilbert  functions. Specifically, if $M$ is a $\mu$-regular bigraded module, then the Hilbert function $H_M(\nu)$ agrees with a polynomial $P_M(\nu)$, termed the Hilbert polynomial of $M$, for all values $\nu \in \left(\mu +\N^2 \right) \setminus \mu$; see \cite[Corollary 2.15.]{MaclaganSmith2}. Furthermore, \cite[Proposition 6.7]{MaclaganSmith} shows that if $I$ is a $B$-saturated ideal defining a finite set of points in \Q, then $\reg(S/I)$ is exactly the set of elements $\mu\in \Z^2$ for which the Hilbert function $H_{S/I}(\mu)$ is equal to the Hilbert polynomial $P_{S/I}(\mu)$.

An important observation from \cite{MaclaganSmith} is that the regularity region of a module $M$ can be estimated from any virtual projective resolution of $M$. We give a version of this result adapted to our setup.

\begin{proposition}[{\cite[Theorem 1.5]{MaclaganSmith}}]
\label{prop:MSbound}
Let  $M$ be a finitely generated bigraded module. If 
$0 \rightarrow  P_3 \rightarrow  P_2 \rightarrow  P_1 \rightarrow  P_0 \rightarrow   M  \rightarrow 0 $
is a virtual projective resolution for $M$ then 
\[
\mathcal{R}=\bigcup_{\sigma:[3]\to[2]}\left( \bigcap_{1\leq i\leq 3} -\be_{\sigma(1)} -\cdots -\be_{\sigma(i)} +\reg(P_i) \right)\subseteq \reg(M).
\]
\end{proposition}

 Unlike the case where the grading group is $\Z$, the minimal free resolution of a bigraded module $M$ does not completely determine its regularity region. This shortcoming is overcome be introducing a related notion of strong regularity developed in \cite{strongregularity}.
\begin{definition}\label{def:strreg}
A bigraded module $M$ is said to be {\em strongly $(a,b)$-regular} if 
\begin{align*}
H^i_{(s,t)}(M)_{(k,k')}=0, & \quad \forall k\geq a+1-i, \forall k'\\
 H^i_{(u,v)}(M)_{(k,k')} =0, & \quad \forall k' \geq b+1-i, \forall k \text{ and }\\
 H^i_{(s,t,u,v)}(M)_{(k+k')} =0, &\quad \forall k+k' \geq a+b+1-i.
 \end{align*}
 We denote by $\reg^s(M)$ the set of all pairs $(a,b)\in \Z^2$ such that $M$ is strongly $(a,b)$-regular.
\end{definition}

It is shown in \cite[Corollary 4.5]{strongregularity} that $\nu\in \reg^s(M)$ implies $\nu \in \reg(M)$. The advantage of strong regularity is that it can be read from the minimal free resolution for the module $M$, Indeed, \cite[Theorem 4.10]{strongregularity} shows that, if for all $i$ the bigraded shifts in the $i$-th homological degree of the minimal free resolution of a module $M$ belong to  
\[
\DReg_i(a,b)=\Z^2_-+\St_{-i}+\mu, \text{ where }Z_- =\{n\in \Z:n\leq0\},
\]
 then $M$ is strongly $\mu$-regular and thus also weakly $\mu$-regular.

\subsection{Eagon-Northcott complex and bigraded regularity}
\label{s:EN}
We follow the notation from the original paper by Eagon and Northcott \cite{ENcomplex}. Let $R$ be a noetherian commutative
ring and let 
\[
\alpha:\bigoplus_{i=1}^r R(-c_i,-d_i)\to \bigoplus_{i=1}^q R(-e_i,-f_i)
\]
 be a bihomogeneous map
where $q,r$ are positive integers with $q\leq r$. Let $I_q(\alpha)$ denote the ideal generated by the maximal minors of any matrix $\phi_\alpha$ representing $\alpha$ with respect to a choice of bases $X_{1},X_{2},\ldots,X_{r}$ for the domain of $\alpha$ and $Y_{1},Y_{2},\ldots,Y_{q}$ for the target of $\alpha$. 
Consider the free graded $R$-modules 
$$K=\bigwedge\left(\bigoplus_{i=1}^r R(-c_i,-d_i)\right)=\bigwedge(X_{1},X_{2},\ldots,X_{r})$$ and 
$$S=\Sym\left(\bigoplus_{i=1}^q R(-e_i,-f_i)\right)=\Sym(Y_{1},Y_{2},\ldots,Y_{q})$$
with $\deg(X_i)=(-c_i,-d_i)$ and $\deg(Y_i)=(e_i,f_i)$ and set $K_i=\bigwedge^i(X_{1},X_{2},\ldots,X_{r})$ and $S_j=\Sym_j(Y_{1},Y_{2},\ldots,Y_{q})$. Let $(-e,-f)=\sum_{i=1}^qq (-e_i, -f_i)$.
The $k$-th row of the matrix $\phi_\alpha=(a_{ij})$ determines a Koszul differential $\Delta_{k}$ on $K$ given by
\[\Delta_{k}(X_{i_{1}}\wedge \cdots \wedge X_{i_{n}})=\sum_{p=1}^{n}(-1)^{p+1}a_{ki_{p}}X_{i_{1}}\wedge \cdots \widehat{X_{i_{p}}}\cdots \wedge X_{i_{n}}.\]
The Eagon-Northcott complex associated to the map $\alpha$, is the complex given by 
\begin{equation*}
0 \rightarrow K_r \otimes_R S_{r-q} \rightarrow \cdots \rightarrow K_{q-i} \otimes_R S_i  \rightarrow K_{q+1}\otimes_R S_1\rightarrow K_q  \rightarrow R(-e,-f) 
\end{equation*}
where the first map  $\bigwedge^q \alpha: K_q\to\bigwedge^q \left( \bigoplus_{i=1}^q R(-e_i,-f_i)\right) = R(-e,-f) $ maps $X_{i_{1}}\wedge \cdots \wedge X_{i_{q}}$ to the maximal minor $\Delta_{i_{1},\ldots, i_{q}}$ of $\phi_\alpha$ determined by the columns $i_{1},\ldots, i_{q}$. The rest of the differentials are specified on the basis elements of $K_{q-i} \otimes_R S_i$ as follows
\[d(X_{i_{1}}\wedge \cdots \wedge X_{i_{q-i}}\otimes Y_{1}^{\nu_{1}}\cdots Y_{q}^{\nu_{q}})=\sum_{j}\Delta_{j}(X_{i_{1}}\wedge \cdots \wedge X_{i_{q-i}})\otimes
Y_{1}^{\nu_{1}}\cdots Y_{j}^{\nu_{j}-1} \cdots Y_{q}^{\nu_{q}}\]
where $\nu_{1}+\ldots+\nu_{s} =i$ and the sum is over those indices $j$ for which $\nu_{j}>0$. With the degree conventions in place this is a complex of free bigraded modules and bidegree $(0,0)$ maps. It is convenient to shift the complex above so that the homological degree 0 component is generated in bidegree $(0,0)$. Henceforth we refer to the shifted version below as the Eagon-Northcott complex $EN(\alpha)$:
\begin{equation}
\label{eq:ENdef}
0 \rightarrow \left(K_r \otimes_R S_{r-q}\right)(e,f)  \rightarrow \cdots \rightarrow \left(K_{q+1}\otimes_R S_1\right)(e,f)   \rightarrow  K_q(e,f)  \rightarrow R.
\end{equation}

The principal application of the Eagon-Northcott complex is in resolving the ideal of minors of matrices $I_q(\phi_\alpha)$ when these ideals have maximum possible height, i.e. $\het \left( I_q(\phi_\alpha)\right) = r-q+1$. The following lemmas are important in establishing the exactness and  computing the homology of the Eagon-Northcott complex in our case of interest. 

\begin{lemma} \label{lem:homology}
Using the notation of \ref{s:EN}, suppose  $r=q+h $. Then
\begin{enumerate}
\item if  $\het(I_q(\alpha))= h$, then $EN(\alpha)$ has  $H_i(EN(\alpha))=0$ for $i\geq 2$.
\item if  $\het(I_q(\alpha))=h+1$, then $EN(\alpha)$ is a resolution for $I_q(\alpha)$.
\end{enumerate}
\end{lemma}
\begin{proof}
By \cite[Theorem~1 Section 5]{ENcomplex}, the homology of the complex $EN(\alpha)$ satisfies
\[\max\{i : H_i\left(EN(\alpha)\right)\neq 0\}=r-q+1-\het(I_q(\alpha))=
\begin{cases}
q+h-q+1-h=1 & \text{ in case } (1)\\
q+h-q+1-h-1=0 & \text{ in case } (2).
\end{cases}
\]
\end{proof}

\begin{remark}\label{rem:subcomplex}
Suppose $q\leq r-1$ and consider a restriction $$\alpha':\bigoplus_{i=1, i\neq i_0}^r R(-c_i,-d_i)\to  \bigoplus_{i=1}^q R(-e_i,-f_i)$$ of the map $\alpha$ defined above, which gives rise to the module $K'=\bigwedge(X_1,\ldots, \widehat{X_{i_0}},\ldots, X_r)$. Since $K'$ is naturally a submodule of $K$,  it follows from $\eqref{eq:ENdef}$ that $EN(\alpha')$ is a subcomplex of $EN(\alpha)$.
In particular, if the degrees of the generators of the free module  $EN(\alpha)_i$ belong to $\DReg_i(a,b)$ then so do the the degrees of the generators of the free module  $EN(\alpha')_i$ implying that
$\reg^s(I_q(\phi_\alpha)) \subseteq \reg^s(I_q(\phi_\alpha'))$. 
In a similar fashion, if  $\mathcal{R}, \mathcal{R}'$ are the
weak regularity regions of $R/I_q(\alpha)$ and $R/I_q(\alpha')$
specified by Proposition~\ref{prop:MSbound}, then $\mathcal{R}\subseteq \mathcal{R}'$.
\end{remark}

\begin{example}
\label{ex:ENci}
We illustrate by showing the  Eagon-Northcott complex when $G$ is a complete intersection. Assume $q=2,r=4$ and $(e_i,f_i)=(k_i,l_i)$ while $(c_1,d_1)=(k_1+k_2,l_1+l_2)$ and $(c_i,d_i)=(a,b)$ for $2\leq i\leq 4$. The bigraded shifts in the Eagon-Northcott complex are illustrated below, based on the degrees of the standard bases of the free modules in the complex \eqref{eq:ENdef}, where $(e,f)=(k_1+k_2,l_1+l_2)$:
\begin{equation*}
\label{ENdegs}
\xymatrix{
0 \ar[r] 
&
\txt{$R(-3a+2k_1,-3b+2l_1)$ \\$\oplus$\\
  $R(-3a+2k_2,-3b+2l_2)$ \\$\oplus$\\
   $R(-3a+k_1+k_2,-3b+l_1+l_2)$}
\ar[r] 
&
\txt{$R(-2a+k_1,-2b+l_1)^{3}$
\\$\oplus$\\$R(-2a+k_2,-2b+l_2)^{3}$
\\$\oplus$\\
$R(-3a+2k_1+k_2,-3b+2l_1+l_2)$
\\$\oplus$\\$R(-3a+k_1+2k_2,-3b+l_1+2l_2)$}
\ar[r]
&
}
\end{equation*}
\begin{equation*}
\xymatrix{
\ar[r]
&
\txt{$R(-a,-b)^{3}$
\\$\oplus$\\
$R(-2a+k_1+k_2,-2b+l_1+l_2)^{3}$}
\ar[r] 
& 
R.
}
\end{equation*}
\end{example}

The following result generalizes Example \ref{ex:ENci}.

\begin{proposition}
\label{prop:ENreg}
Let $\alpha:\bigoplus_{i=1}^{n-1} R(-c_i,-d_i)\oplus R(-a,-b)^3\to \bigoplus_{j=1}^{n}R(-e_j,-f_j)$ be a bidegree preserving map and set $(c,d)=\sum_{i=1}^{n-1} (c_i,d_i)$ and $(e,f)=\sum_{j=1}^{n}(e_j,f_j)$. Then the degrees of the minimal generators for the free $R$-modules in the complex $EN(\alpha)$, listed by homological degree, are as follows
\medskip

\begin{tabular}{|c|c|}
\hline
\text{degree} & \text{shifts} \\
\hline
0 &  (0,0) \\
1 &  $(a+c-e, b+d-f),(2a+c-e-c_i,2b+d-f-d_i)$\\
& $(3a+c-e-c_i-c_j,3b+d-f-d_i-d_j), i\neq j$\\
2 & $(2a+c-e-e_j,2b+d-f-f_j),(3a+c-e-c_i-e_j,2b+d-f-d_i-f_j)$\\
3 &  $(3a+c-e-e_i-e_j,3b+d-f-f_i-f_j)$\\
\hline
\end{tabular}
\medskip

\noindent In particular, if $a\geq e_j$ for some $1\leq j\leq n$ and $b\geq f_j$ for some $1\leq j\leq n$ and $EN(\alpha)$ is a virtual resolution for a module $R/I_q(\phi_\alpha))$ then the bigraded regularity of $R/I_q(\phi_\alpha))$ can be estimated by
\[
\mathcal{R}(\alpha)=\left(3a+c-e-\!\!\min_{1\leq i\leq j\leq n}(e_i+e_j),3b+d-f-\!\!\min_{1\leq i\leq j\leq n}(f_i+f_j)\right)+\St_{-3}+\N^2\subseteq \reg\left(R/I_q(\phi_\alpha))\right).
\]

\end{proposition}
\begin{proof}
The  shifts listed in the table follow from the graded structure of the complex \eqref{eq:ENdef}. 

Denoting by $P_i$ the free module in the $i$-th homological degree in $EN(\alpha)$ 
we claim that $\reg(P_{i+1})\subseteq\reg(P_i)+(1,1)$ for $0\leq i\leq 4$. Using the fact that for any two modules $U,V$ $\reg(U\oplus V)=\reg(U)\cap\reg(V)$ (\cite[Lemma 7.1]{MaclaganSmith}) one can easily compute the regularity of a graded free $R$-module $\bigoplus_{i=1}^q R(-m_i,-n_i)=\left(\max\limits_{1\leq i\leq q}m_i,\max\limits_{1\leq i\leq q}n_i\right)+\N^2$.
Thus, to establish the claim it is sufficient to show that the maximum of the first components of the degrees listed in row $i$ of the table above is strictly smaller than the maximum of the first components of the degrees listed in row $i+1$ of the table and the analogous statement for the second components. For $i=0$ this is clear, so we assume $i>0$. Notice that $a\geq e_j$ for some $1\leq j\leq n$ and $b\geq f_j$ for some $1\leq j\leq n$ ensures that each component of the degrees listed in row $i$ of the table above is strictly smaller than some component of the degrees listed in row $i+1$, which establishes the claim.

In particular,  $\reg(P_{i+1})\subseteq\reg(P_i)+(1,1)$ implies that 
\[
-\be_{\sigma(1)} -\cdots -\be_{\sigma(i)} +\reg(P_i)
\subseteq 
-\be_{\sigma(1)} -\cdots -\be_{\sigma(i)} -\be_{\sigma(i+1)} +\reg(P_{i+1}) \text{ for } 0\leq i\leq 2.
\]
The statement of Proposition \ref{prop:MSbound} can now be simplified to say
\[
\bigcup_{\sigma:[3]\to[2]}\left( -\be_{\sigma(1)} -\be_{\sigma(2)} -\be_{\sigma(3)}+\reg(P_3) \right)\subseteq \reg(M), \text{ i.e. }
 \reg(P_3)+\St_{-3} \subseteq \reg(M).
\]
Using the explicit formula for the regularity of a free module deduced above yields the desired estimate
\[
\left(3a+c-e-\min_{1\leq i\leq j\leq n}(e_i+e_j),3b+d-f-\min_{1\leq i\leq j\leq n}(f_i+f_j)\right)+\St_{-3}+\N^2\subseteq \reg\left(R/I_q(\phi_\alpha))\right).
\]
\end{proof}

\begin{remark}
\label{rk:ENregspecialized}
All the results of this section continue to hold verbatim for $T$-modules. In particular, if $\alpha:\bigoplus_{i=1}^{n-1} T(-c_i,-d_i)\oplus T(-a,-b)^3\to \bigoplus_{j=1}^{n}T(-e_j,-f_j)$ is a bidegree preserving map, $(c,d)=\sum_{i=1}^{n-1} (c_i,d_i)$ and $(e,f)=\sum_{j=1}^{n}(e_j,f_j)$, then the region $\mathcal{R}(\alpha)$ of Proposition~\ref{prop:ENreg} is contained in the regularity region of $I_q(T/I_q(\phi_\alpha))$, provided that the Eagon-Northcott complex is a virtual projective resolution for this module.

Note that the region $\mathcal{R}(\alpha)$  only depends on the numerical information regarding the degrees in which the domain and target of the map $\alpha$ are generated and not on the rule defining $\alpha$. In particular applying an evaluation map to the source  and target of $\alpha$ induces an $R$-linear map $e_c(\alpha):\bigoplus_{i=1}^{n-1} R(-c_i,-d_i)\oplus R(-a,-b)^3\to \bigoplus_{j=1}^{n}R(-e_j,-f_j)$ such that $\mathcal{R}(e_C(\alpha))=\mathcal{R}(\alpha)$.
\end{remark}

\section{\large Effective computation of the residual resultant}\label{s:compresultant}
\label{s:resultantcomputation}

\subsection{Virtual resolutions for effective computations}

Let $G\subseteq R$ be a bihomogeneous ideal.  For $0\leq i\leq m$, let $(a_i,b_i)\in \N^2$  and set $C=k[C_{ij}^{\alpha}: 0\leq i\leq m,1\leq j\leq n] $ where for each pair $i,j$, the index $\alpha$ enumerates the elements  $m_{\alpha}$ of a monomial basis of  $R_{(a_i-k_j,b_i-l_j)} $. Define $H_{ji}=\sum_{\alpha} C_{ij}^{\alpha}m_{\alpha}$, $F_i=\sum_{j=1}^n H_{ji}g_j$, so  $H_{ji}\in T_{(a_i-k_j,b_i-l_j)}$ and $F_i\in T_{(a_i,b_i)}$. This can be written concisely as 
\begin{equation}
\label{eq:FandPsi}
\begin{bmatrix} F_0 & \cdots & F_m \end{bmatrix} =  \begin{bmatrix} g_1&\ldots & g_n \end{bmatrix} \Psi , \text{ where } \Psi=[H_{ji}]_{1\leq j\leq n,\\0\leq i\leq m}\in \cM_{n\times m}(T).
\end{equation}
Lastly, set $\F=(F_0,\ldots, F_m)$ and notice that the previous equation gives the containment $\F\subseteq G$. We study the algebraic counterpart of the residual resultant developed in Section \ref{s:resultant}. 
 As mentioned previously, we denote by $h_{ji}$ and $\psi$ the images of $H_{ji}$ and $\Psi$ under any evaluation homomorphism $e_c:T\to R$. 

We aim  to express the residual resultant for the pair of ideals $\F,G$ in terms of the minimal free resolution for the residual ideal $I=\F:_TG$. In turn, we will approximate this resolution by a virtual projective resolution  based on the structure matrix $\Psi$ defined above as well as the syzygy matrix $\varphi$ for $G$. The  syzygy matrix for $G$ is determined up to change of basis by a minimal presentation $R^\ell \stackrel{\varphi}{\longrightarrow} R^n \to G \to 0$. We start with a lemma that relates the ideal $\F:_TG$ to the ideal of maximal minors of the matrix $\varphi\oplus \Psi\in \cM_{n\times(\ell+m)}(T)$.

\begin{lemma}\label{lem:phipsi}
Let $\F=(F_0,\ldots,F_m)\subseteq G=(g_1,\ldots,g_n)$ be homogeneous ideals in $T$  with the sets of generators 
of the two ideals related by
\[\begin{bmatrix} F_0 & \cdots & F_m \end{bmatrix} =  \begin{bmatrix} g_1&\ldots & g_n \end{bmatrix} \Psi .\]
Let $\varphi$ denote the $n\times \ell$ matrix of syzygies of $G$. 
Then  the following hold

\begin{enumerate}
\item  $I_n(\varphi\oplus \Psi)\subseteq \Ann(\coker(\varphi\oplus \Psi)) =\F:_TG.$
\item if  $\h(\F:_TG)=m-n+\ell+2$, then equality holds in the above containment.  
\end{enumerate}
\end{lemma}
\begin{proof}
Computing ranks along the exact sequence $R^\ell\stackrel{\varphi}{\longrightarrow} R^n \to R\to R/G\to 0$ gives $\ell\geq n-1$, thus $\ell+m\geq n$, hence $I_n(\varphi\oplus \Psi)$ is the ideal of maximal minors of $\varphi\oplus \Psi$.
Note that $\varphi\oplus \Psi$ appears in the following bigraded presentation for  $G/\F$:
\[
\bigoplus_{i=0}^{m+\ell+1} T(-e_i,-d_i) \stackrel{\varphi\oplus \Psi}{\longrightarrow} \bigoplus_{i=1}^{n} T(-k_i,-l_i) \to G/\F \to 0.
\]
 A theorem of Buchsbaum-Eisenbud \cite{BE} on Fitting ideals, applied to the presentation above,  gives the containment below,  with equality instead of the rightmost containment when $\het\left(I_n(\varphi\oplus \Psi) \right)=m-n+\ell+2$:
  \begin{equation}\label{BE}
  \Ann\left(\coker(\varphi\oplus \Psi)\right)^n\subseteq I_n(\varphi\oplus \Psi) \subseteq \Ann\left(\coker(\varphi\oplus \Psi)\right).
  \end{equation}
  
Combining the containment above and the identity 
\[
\Ann\left(\coker(\varphi\oplus \Psi)\right)=\Ann (G/\F) = \left(\F :_T G\right)
\]
 gives the first statement of the lemma.
Furthermore, if $\het(\F:_TG)=m-n+\ell+2$, 
 the containment  (\ref{BE}) and the generalized principal ideal theorem (see \cite[Exercise 10.9]{Eisenbud}) $\het\left(I_n(\varphi\oplus \Psi) \right)\leq m-n+\ell+2$ , yield   $\het\left(I_n(\varphi\oplus \Psi) \right)= m-n+\ell+2$, which gives the second statement of the lemma. 
\end{proof}

Note that the identity $I_n(\varphi\oplus \Psi)=\F:_TG$ can hold  even if the hypothesis of statement (2) above is not met, as illustrated in Example~\ref{ex:2}.

\begin{cor}\label{cor:phipsi}
The statement of the lemma holds over the ring $R$ whenever $F_0,\ldots, F_m$ and $\Psi$ are specialized  via evaluation to $R$.
\end{cor}

We exploit the close relation between $I_n(\varphi\oplus \Psi)$ and $\F:_TG$ established in Lemma~\ref{lem:phipsi} to obtain a virtual resolution of $\F:_TG$. First, due to Proposition~\ref{prop:HBpoints} we may assume that $G$ is an ideal with a Hilbert-Burch resolution provided the degrees of the generators of $\F$ are high enough. The exact meaning of this reduction is made precise in the following proposition.

\begin{lemma}
\label{lem:GwithHBres}
Suppose that $G$ defines a zero-dimensional subscheme of $\Q$ and $\F$ is an arbitrary ideal of $T$. Then there exists an ideal $G'$ of $R$ that has the following properties:
\begin{enumerate}
\item $V(G)=V(G')$,
\item $G'$ has a Hilbert-Burch resolution,
\item $\mathrm{Res}_{\mathcal{G},\{(a_{i},b_{i})\}_{i=0}^2}=\mathrm{Res}_{\mathcal{G}',\{(a_{i},b_{i})\}_{i=0}^2}$ for $(a_i,b_i)$ satisfying the condition in Proposition \ref{prop:tpres},
\item $(\F:G)^{\rm sat}=(\F:G')^{\rm sat}$ where saturation is taken with respect to the ideal $\B$ of $T$,
\item a complex $\mathbf{F}$ of free $T$ modules is a  virtual projective resolution for $\F:G$ if and only if $\mathbf{F}$ is  a virtual projective resolution for $\F:G'$ as well.
\end{enumerate}
Moreover, if the reduced subscheme of $\Q$ defined by $G$ consists of $r$ general points and $\F\subseteq G$ is an ideal of $T$ such that the generators of $\F$ have bidegrees lying in the interior of the region $\mathcal{S}$ shown in Figure~\ref{fig:general} then $\F\subseteq G'$.
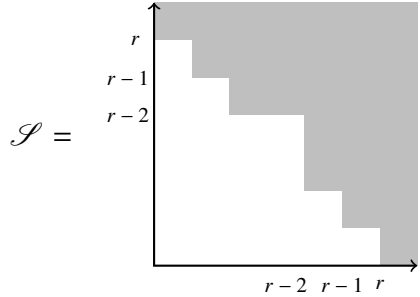
\begin{figure}[h!]
\label{fig:general}
\begin{tikzpicture}[scale=0.5]
\node  at (-3,3.5){$\mathcal{S}=$};
   \path [fill=lightgray] (0,6) rectangle (7,7);
   \path [fill=lightgray] (1,5) rectangle (7,7);
   \path [fill=lightgray] (2,4) rectangle (7,7);
   \path [fill=lightgray] (6,0) rectangle (7,7);
    \path [fill=lightgray] (5,1) rectangle (7,7);
    \path [fill=lightgray] (4,2) rectangle (7,7);
   \draw [<->] [thick] (0,7) -- (0,0) -- (7,0);
   \node at (-0.5,6) {\tiny$r$};
     \node at (-0.7,5) {\tiny$r-1$};
   \node at (-0.7,4) {\tiny$r-2$};
 \node at (6, -0.5) {\tiny$r$};
     \node at (5, -0.5) {\tiny$r-1$};
   \node at (3.5,-0.5) {\tiny $r-2$};
  \end{tikzpicture}
   \vspace{-0.1cm}
\caption{Region $\mathcal{S}$ referred to in Lemma~\ref{lem:GwithHBres}}
 \vspace{-0.2cm}
\end{figure}
\end{lemma}
\begin{proof}
Let $G'$ be the ideal given by Corollary~\ref{cor:HBpoints}, which establishes that it satisfies properties (1) and (2) listed above. Note that property (1) is equivalent to $G^{\rm sat}={G'}^{\rm sat}$ and therefore $\mathcal{G}=\mathcal{G}'$, which yields property (3) tautologically.

For (4), consider $\fp\in \Spec(T)$. If $B\not \subseteq \fp$ (equivalently $\B\not \subseteq \fp$) then the equality $G^{\rm sat}={G'}^{\rm sat}$ implies that $G_\fp=G'_\fp$ and therefore we have 
$(\F:G)_\fp=\F_\fp:G_\fp=\F_\fp:G'_\fp=(\F:G')_\fp$, which is equivalent to $(\F:G)^{\rm sat}=(\F:G')^{\rm sat}$.
For (5), recall that $\mathbf{F}$ is a virtual resolution of $\F:G$ if and only if $\left(H_0(\mathbf{F})\right)^{\rm sat}=(\F:G)^{\rm sat}=(\F:G')^{\rm sat}$ and $H_i(\mathbf{F})$ is $\B$-torsion for $i>0$.

When the reduced locus of $G$ consists of $r$ general points then $G'$ can be taken to have one of the two types of resolutions presented in Example \ref{ex:generalpts} or those obtained from the ones presented by interchanging the two coordinates of each bidegree. By Proposition~\ref{prop:MSbound} the region $\mathcal{S}$ is contained in the union of the regularity regions of the two possible cyclic modules $R/G'$ afforded by the value of $r$.  Note that $\mathcal{S}$ is also contained in the regularity region of $R/G$ because the resolutions in Example \ref{ex:generalpts} are virtual projective resolutions for $R/G$. Because $G_\nu=G^{\rm sat}_\nu={G'}^{\rm sat}_\nu=G'_\nu$ for any $\nu=\mu+(i,j)$ with $\mu\in\mathcal{S}$ and $i,j\in \N$ with $i+j>0$ (denote this by $\nu\in\mathcal{S}^0$) we have that $F\subseteq \bigcup_{\nu\in\mathcal{S}^0} G_\nu$  implies $F\subseteq \bigcup_{\nu\in\mathcal{S}^0} G'_\nu \subseteq G'$.
\end{proof}

The existence of a Hilbert-Burch resolution for $G$ is a key ingredient in our results and from this point on we assume that $G$ satisfies this property. We further assume that $m=2$ since this is the setup for a residual resultant over $\Q$. Under these conditions $\varphi\in M_{n\times(n-1)}(T)$ and the matrix $\varphi\oplus \Psi $ in Lemma~\ref{lem:phipsi} is a $n\times(n+2)$ matrix. 

\begin{proposition} \label{prop:virtualT}
Assume that $G\subseteq R$ has a Hilbert-Burch resolution,  $\F=(F_0, F_1, F_2)$, and suppose that for every $\fp\in \Spec(R)\setminus \cB$ with $\h (\fp)=2 $ there is an equality $F_{\fp}=G_{\fp}$. Then the Eagon-Northcott complex $EN(\varphi \oplus \Psi)$  is a virtual resolution for  the module $T/I_n(\varphi\oplus \psi)$.
\end{proposition}
\begin{proof}
Throughout this proof, let $\Min(I)$ denote the set of minimal primes of an ideal $I$. 

Recall from Remark~\ref{rem:incidence} that the incidence variety $W\subset \left( Q\setminus Z \right)\times \prod_{i=0}^{2}V_i$ has codimension three. Since $V(F:_TG)\subseteq W\cup Z$, 
it follows that there is a containment 
\[
\Min_T(F:_TG)\subseteq \Min_T(I_W) \cup \cB_T\cup \Min_T(G),
\]
where $\cB_T$ and  $\Min_T(G)$ are the set of primes in $\cB$ and $ \Min_R(G)$ respectively extended to $T$. Since $F_{\fp}=G_{\fp}$ holds for any $\fp\in \Spec(R)\setminus \cB$ with height of $\fp$ equal to two, it follows that in fact any prime of $ \Min_T(G)$ that is also an associated prime of $F:_TG$ is in $\cB_T$, thus the containment above reduces to 
\begin{equation}
\label{eq:mincolon}
\Min_T(F:_TG)\subseteq \Min_T(I_W) \cup \cB_T.
\end{equation}

Consider the $T$-module $H=\bigoplus_{i=0}^{n+2} H_i\left(EN(\varphi\oplus \Psi)\right)$.
From \cite[Theorem A2.59]{Eisenbud} it follows that $I_n(\varphi\oplus \psi)$ annihilates $H$, therefore there is a containment 
\begin{equation}
\label{eq:minH}
\Min_T(H)\subseteq \Min_T\left(I_n\left(\varphi\oplus \Psi\right)\right).
\end{equation}
The containments $(F:_TG)^n\subset I_n\left(\varphi\oplus \Psi\right)\subseteq (F:_TG)$ noted in the proof of Lemma~\ref{lem:phipsi} imply $\sqrt{F:_TG}=\sqrt{I_n\left(\varphi\oplus \Psi\right)}$ and hence $\Min_T\left(I_n\left(\varphi\oplus \Psi\right)\right)=\Min_T(F:_TG)$. Therefore, from equations \eqref{eq:mincolon}, \eqref{eq:minH} we deduce
\begin{equation}
\label{eq:combined}
\Min_T(H)\subseteq \Min_T\left(I_n\left(\varphi\oplus \Psi\right)\right)=\Min_T(F:_TG)\subseteq \Min_T(I_W) \cup \cB_T.
\end{equation}
Let $P$ be any ideal of $T$ of height at least 3; in particular this applies to any $P\in \Min_T(I_W)$ since the codimension of $W$ is 3 by Remark \ref{rem:incidence}. Then the complex $EN\left(\varphi\oplus \Psi\right)\otimes_T T_P=EN_{T_P}\left(\varphi\oplus \Psi\right)$ is exact by Lemma \ref{lem:homology} (2) because $\het\left(I_n\left(\varphi\oplus \Psi\right)_P\right)=\het(P)\geq 3$. It follows that $H_P=0$ and therefore $P$ is not in the support of $H$, so $P\not\in \Ass_T(H)$. This shows that the associated primes of $H$ have height 2 and further reduces equation \eqref{eq:combined}  to 
\begin{equation}
\label{eq:simple}
\Min_T(H)\subseteq \Min_T\left(I_n\left(\varphi\oplus \Psi\right)\right)\subseteq\cB_T.
\end{equation}
Therefore $\Ann_T(H)= Q_1 \cap Q_2$ where
 $Q_1$ is $\langle s,t\rangle$-primary and $Q_2$ is $\langle u,v\rangle$-primary.
 Now $\langle s,t\rangle^a\subseteq Q_1$ for some $a\geq 0$ and similarly $\langle u,v\rangle^b\subseteq Q_2$ for some $a\geq 0$, hence for $m\geq\max\{a,b\}$ we have the desired conclusion
\[
B^m=\langle s,t\rangle^m\cap \langle u,v\rangle^m\subseteq Q_1\cap Q_2\subseteq \Ann_T(H).
\]
\end{proof}

\begin{remark}
Example \ref{ex:2} illustrates the fact that it is possible for the Eagon-Northcott complex in Corollary \ref{cor:virtualR} to be a virtual projective resolution while not being a resolution, i.e. not being exact.
\end{remark}

\begin{remark}
In the setup of this section, where $F_i=\sum_{j=1}^n H_{ji}g_j$ for $0\leq i\leq m$, $H_{ji}=\sum_{\alpha} C_{ij}^{\alpha}m_{\alpha}$ and  $\alpha$ runs over the elements  $m_{\alpha}$ of a monomial basis of  $R_{(a_i-k_j,b_i-l_j)} $, the hypothesis that there is an equality $\F_{\fp}=G_{\fp}$  for every $\fp\in \Spec(R)\setminus \cB$ with $\h (\fp)=2 $ holds true whenever $G$ is locally a complete intersection (see Lemma~\ref{lem:hypothesisgeneric}). However, we prefer to state Proposition~\ref{prop:virtualT} including this hypothesis, since we shall use it in a slightly more general context in section~\ref{s:implicitization} and also to draw a closer analogy with the following corollary.
\end{remark}

\begin{cor}\label{cor:virtualR}
Suppose that $G$ has a Hilbert-Burch resolution and the ideal $F=(f_0,f_1,f_2)$ arising by specializing the coefficients of $F_0,F_1,F_2$ to values in $k$ satisfies $F^{\rm sat}=G^{\rm sat}$. Denote by $\psi$ the corresponding specialization of the matrix $\Psi$ in the setup at the beginning of this section. Then the Eagon-Northcott complex $EN(\varphi \oplus \psi)$ over $R$ is a virtual resolution for  the module $R/I_n(\varphi\oplus \psi)$.
\end{cor}
\begin{proof}
 By Corollary~\ref{cor:phipsi}, the conclusion of Lemma~\ref{lem:phipsi} still holds for $f_0,f_1,f_2$. The hypothesis $F^{\rm sat}=G^{\rm sat}$ implies that $F_\fp=G_\fp$ for all $\fp\in \Spec(R)\setminus\cB$ and thus $\Min(F:_RG)\subseteq \cB$. Therefore the proof of Proposition~\ref{prop:virtualT} starting at equation \eqref{eq:simple} applies to show that $EN(\varphi\oplus\psi)$  is $B$-torsion as a complex over $R$.
\end{proof}

\begin{lemma}
\label{lem:hypothesisgeneric}
Assume that $G$ is a locally complete intersection ideal and $F_i=\sum_{j=1}^n H_{ji}g_j$ for $0\leq j\leq m$, where $H_{ji}=\sum_{\alpha} C_{ij}^{\alpha}m_j^{\alpha}$ and $m_j^{\alpha}$ runs over the elements of a monomial basis of  $R_{(a-k_j,b-l_j)}$ for some $(a,b)\in \N^2$. Then there is an equality $\F_{\fp}=G_{\fp}$  for every $\fp\in \Spec(R)\setminus \cB$ with $\h (\fp)=2$.
\end{lemma}
\begin{proof}
Let $\fp\in \Spec(R)\setminus \cB$ be an ideal with $\h (\fp)=2$. We show that for any pair $i,j$ we have $H_{ji}\not\in \fp$. Assume the contrary, fix $\alpha_0$ in the indexing set of monomials in $R_{(a-k_j,b-l_j)}$ and consider the prime ideal $\fq=\fp+(C_{ij}^\alpha: \alpha\neq \alpha_0)$. Then $H_{ji}\in \fp$ implies that $C_{ij}^{\alpha_0}m_j^{\alpha_0}\in \fq $ and since $C_{ij}^{\alpha_0}\not \in \fq$ this yields $m_j^{\alpha_0}\in \fq$, which in turn implies that $m_j^{\alpha_0}\in \fq\cap R=\fp$ for any $\alpha_0$. We deduce that 
$R_{(a-k_j,b-l_j)}=\langle s,t \rangle ^{a-k_j} \cap \langle u,v \rangle^{b-l_j}\subseteq \fp$ and consequently $\fp \in \cB$, a contradiction. Therefore the elements $H_{ij}$ become units in $T_\fp$. 
Since $G_\fp$ is a complete intersection with $\dim_{k(\fp)}G_\fp/G^2_\fp=2$ and $\F_\fp\subseteq G_\fp$ is generated by 3 elements which are pairwise independent in  $G_{\fp}/G^2_\fp$, the equality $\F_\fp=G_\fp$ follows. 
\end{proof}

\subsection{A matrix representation for the residual resultant}\label{s:matrixrep}
The computation of the residual resultant hinges on the following proposition, which identifies a matrix whose rank drops when evaluated at any point of the residual resultant. In an alternate terminology, the following proposition gives a matrix representation for the residual resultant.
\begin{proposition}
\label{prop:surjective}
Let $g_1, \ldots, g_n$ and $f_0, f_1, f_2$ be polynomials in $R$  with $f_i\in R_{(a_i,b_i)}$  related by the identities $f_i=\sum_{j=1}^n h_{ji}g_j$ . Set $G=(g_1\ldots,g_n)$,  $\psi=[h_{ji}]$, and assume that $G$ has a Hilbert-Burch syzygy matrix $\varphi$. Let $\theta$ be a presentation map for the cyclic module $R/I_n(\varphi\oplus \psi)$. The following statements are equivalent
\begin{enumerate}
\item $\Res_{\mathcal{G},\{(a_i,b_i)\}_{i=0}^2}(f_0,f_1,f_2)\neq 0$,
\item $V (I_n(\varphi\oplus \psi)) = \emptyset$,
\item the restriction of the map $\theta$ to degree $\nu$  is surjective for all degrees $ \nu = \mu+(p,p')$  such that $\mu\in \cR(\varphi\oplus\psi), (p,p')\in\N^2$ and $p+p'>0$.
\end{enumerate}
\end{proposition}
\begin{proof}
$(1)\Leftrightarrow (2):$
By Proposition \ref{prop:tpres}, the condition $\Res_{G,\{(a_i,b_i)\}_{i=1}^s}(f_0,f_1,f_2)\neq0$  is equivalent to 
$F^{\rm sat}= G^{\rm sat}$, which is equivalent to $F^{\rm sat} :_R G^{\rm sat}= R$. In view of 
Corollary~\ref{cor:phipsi},  this translates to $\sqrt{I_n(\varphi\oplus\psi)^{\rm sat}}=\sqrt{(F:_RG)^{\rm sat}}=\sqrt{(F^{\rm sat}:_RG^{\rm sat})}= R$ that is, $V (I_n(\varphi\oplus \psi)) = \emptyset$. 


$(1)\Rightarrow (3):$ By Proposition~\ref{prop:tpres} $\Res_{G,\{(a_i,b_i)\}_{i=0}^2}(f_0,f_1,f_2)\neq 0$ implies that $F^{\rm sat}=G^{\rm sat}$, whence Corollary~\ref{cor:virtualR} implies that the Eagon-Northcott complex is a virtual projective resolution for $R/I_n(\varphi\oplus\psi)$ and this module is $\mu$-regular for $\mu\in \cR(\varphi\oplus\psi)$.  Since $V(I_n(\varphi\oplus\psi))=V(F:_RG)=\emptyset$ by hypothesis and $(1)\Rightarrow (2)$, we deduce from \cite[Corollary 2.15]{MaclaganSmith2} that $H_{R/I_n(\varphi\oplus\psi)}(\nu)=0$ for bidegrees $\nu=\mu+(p,p')$  such that  $(p,p')\in\N^2$ and $p+p'>0$. 
Since the cokernel of the restriction of the map $\theta$ to degree $\nu$ is $R/I_n(\varphi\oplus\psi)_\nu$, and by the previous considerations $R/I_n(\varphi\oplus\psi)_\nu=0$ we deduce that this map is surjective.

For $(3)\Rightarrow (2)$ we prove the contrapositive. Suppose that $V(I_n(\varphi\oplus\psi)$ is not empty. Due to the equality $V(I_n(\varphi\oplus\psi))=V(F:G)$, there
exists a point  $\xi\in V(F)\setminus V(G)$. Evaluating the following identity encompassing the expressions $f_i=\sum_{j=1}^n h_{ji}g_j$ and the fact that $\varphi$ is a syzygy matrix for $G$ at $\xi$
$$\begin{bmatrix} 0 & \cdots & 0 &f_0 & f_1& f_2 \end{bmatrix} = \begin{bmatrix} g_1& \dots & g_n\end{bmatrix} \begin{bmatrix} \varphi \oplus \psi \end{bmatrix},$$
shows that the rank of the  matrix  $\varphi \oplus \psi$ evaluated at $\xi$  is not maximal ($<n$). Hence all the maximal minors of $\varphi \oplus \psi$ vanish at $\xi$. Since all these minors generate $I_n(\varphi \oplus \psi)$ we deduce that for arbitrary $\nu\in \N^2$ any polynomial in the image of the map $\theta_\nu$  vanishes at $\xi$. Since for any point $\xi$ there exist polynomials in $R_\nu$ that do not vanish at $\xi$, it follows that the map $\theta_\nu$ is not surjective.
\end{proof}

\begin{remark}\label{rem:regularity}
Proposition~\ref{prop:surjective} relates the nonvanishing of  $\Res_{\mathcal{G},\{(a_i,b_i)\}_{i=0}^2}(f_0,f_1,f_2)$ to
the presentation of the module $R/I_n(\varphi\oplus \psi)$ restricted to any bidegree in the interior of the region $\cR(\varphi\oplus\psi)$ described in Proposition~\ref{prop:ENreg}. Note that by Remark~\ref{rk:ENregspecialized} this region is stable under specialization, that is $\cR(\varphi\oplus\psi)=\cR(\varphi\oplus\Psi)$.

We now proceed to convert Proposition~\ref{prop:surjective} into an effective computational tool.  
\end{remark}

In order to make the matrix representation for the residual resultant explicit we recall the first map of the Eagon-Northcott complex \eqref{eq:ENdef} associated to the matrix $\varphi\oplus\Psi$ over $T$, 
\[
d=\bigwedge^q\left(\varphi\oplus\Psi\right):\bigoplus_{\{i_1,\ldots,i_q\} \subset [r]  }T X_{i_1}\wedge\cdots \wedge X_{i_q}\to T, X_{i_1}\wedge\cdots \wedge X_{i_q} \mapsto \Delta_{i_1,\ldots,i_q}. \]
Here $\Delta_{i_1,\ldots,i_q}$ is the maximal minor of $\varphi\oplus\Psi$ corresponding to the columns $i_1,\ldots,i_q$ and the $T$-module generated by $X_{i_1}\wedge\cdots \wedge X_{i_q}$ is generated in degree $\deg(\Delta_{i_1,\ldots,i_q})$. For $\nu\in \Z^2$, let $d_\nu$ denote the map $d$ restricted
to bidegree $\nu$. Since for any bidegree $\nu$, $T_{\nu}=C\otimes_k R_{\nu}$ is a free $C$-module, we obtain a map of finitely generated free $C$-modules.
\[ 
d_\nu:\bigoplus_{\{i_1,\ldots,i_q\} \subset [r]  }T_{\nu-\deg(\Delta_{i_1,\ldots,i_q})}\to T_\nu.
\]
An explicit matrix representing the map $d_\nu$ can be obtained in four  steps:
\begin{enumerate} 
\label{steps}
\item  fix a basis for the  vector space $\bigoplus_{\{i_1,\ldots,i_q\} \subset [r]  }T_{\nu-\deg(\Delta_{i_1,\ldots,i_q})}$, 
\item apply the map $d_\nu$ to this basis,
\item fix a basis for $R_\nu$ and express the result of step (2) in terms of this basis as vectors with entries in $C$,
\item form a matrix with entries in $C$ denoted $\Theta_\nu$ having these vectors as columns. 
\end{enumerate}
Note that for step one, a standard basis of this vector space consists of elements $mX_{i_1}\wedge\cdots \wedge X_{i_q}$ such that $m$ is a monomial in $R$ with $\deg (m)=\nu-\deg(\Delta_{i_1,\ldots,i_q})$ for some $\{i_1,\ldots,i_q\}\subseteq [n]$. Then in step (2) one obtains $d_\nu(m)=m\cdot \Delta_{i_1,\ldots,i_q}$. 

For any bidegree $\nu\in\Z^2$ we denote by $\theta_\nu$ the image of the matrix $\Theta_\nu$ defined above under an evaluation homomorphism. According to part (3) of Proposition~\ref{prop:surjective}, from  this point onward we let $\nu$ be a bidegree such that $\nu= \mu+(p,p')$ with $\mu\in \mathcal{R}(\varphi\oplus\Psi)$ and $(p,p')\in \N^2, p+p'>0$. When this holds we say that $\nu$ is in the interior of $\mathcal{R}(\varphi\oplus\Psi)$.

\begin{proposition}\label{pro:multiple}
If  $\nu$ is in the interior of $\mathcal{R}(\varphi\oplus\Psi)$, then any nonzero minor of size $\dim_{k}(R_\nu)$ of the matrix $\Theta_\nu$ is a multihomogeneous polynomial in the coefficients $C_{ij}^{\alpha}$ of $F_0,F_1,F_2$ and a multiple of $\Res_{\mathcal{G},\{(a_i,b_i)\}_{i=0}^2}$.
\end{proposition}

 In light of Proposition~\ref{prop:surjective} this proof follows along the lines of the argument in \cite{residualBuse}. We include the details for completeness.
 
\begin{proof}
First observe that any minor $\rho$ of the matrix $\Theta_{\nu}$ is multihomogeneous  in the coefficients of each $F_i$ for $i=0,1,2$. Indeed, if $F_i$ is multiplied by a scalar $\lambda\in k$ then the same is true for the column in $\varphi\oplus \Psi$ that corresponds to the coefficients of $\lambda F_i$. Consequently any column in $\Theta_{\nu}$ containing the coefficients of an element  $m\cdot \Delta_{i_1,\ldots,i_s}$ such that
 $\Delta_{i_1,\ldots,i_s}$ involves a column corresponding to the coefficients of $\lambda F_i$  is multiplied by a factor of $\lambda^n$. This implies that $\rho$ is homogeneous of degree $n\cdot \dim_{k}R_{\nu}\cdot d$ in the coefficients of $F_i$, where $d$ is the number of columns that appear in the the submatrix of $\Theta_{\nu}$ that have a factor of $\lambda^n$.

Next, fix $\rho$ to be a maximal minor of $\Theta_{\nu}$. We want to show that $\rho$ vanishes at every point where the resultant vanishes, for this implies $\rho$ is a multiple of $\Res_{\mathcal{G},\{(a_i,b_i)\}_{i=0}^2}$. 
Let $Q=\Q$ and let $\widetilde{Q}$ be the blow-up of $\Q$ along the sheaf of
ideals associated to $G$. Define $\widetilde{Q^0}=\widetilde{Q}\setminus E$ where $E$ is the exceptional divisor in $\widetilde{Q}$. Let 
\[Z^0=V\left(\Res_{\mathcal{G},\{(a_i,b_i)\}_{i=0}^2}\right)=\{c=(c_{ij}): \exists x\in \widetilde{Q^0}, \pi^{\ast}(f_0)=
\pi^{\ast}(f_1)=\pi^{\ast}(f_2)=0\},
\]
 i.e $Z^0$ is the set of coefficients such that the pullbacks of the sections $f_0,f_1,f_2$ have a common root outside the exceptional divisor $E$. Suppose there is a choice of coefficients $c\in Z^0$ such that $e_c(\rho)\neq 0$. This implies that $\theta_\nu$ is surjective because $\rho$ is
a maximal nonvanishing minor of size $\dim_{k}(R_\nu)$. However, since $c\in Z^0$, the specialized sections $f_0,f_1,f_2$   have a common root in $V(F)\setminus V(G) $ by Proposition~\ref{prop:tpres}. Using the equivalence $(1)\ifi(3)$ of Proposition~\ref{prop:surjective} this implies that $\theta_\nu$ cannot be surjective, a contradiction. Therefore $e_c(\rho)=0$ and since $c\in Z^0$ was arbitrary $\rho$ vanishes on $Z^0$. As $\widetilde{Q^0}$ is dense in $\widetilde{Q}$, $Z^0$ is also dense in $Z=\{c=(c_{ij}): \exists x\in \widetilde{Q}, \pi^{\ast}(f_0)= \pi^{\ast}(f_1)=\pi^{\ast}(f_2)=0\}$. Consequently, $\rho$ vanishes on $Z$, i.e. $\rho$ vanishes at all the points where $\Res_{G,\{(a_i,b_i)\}_{i=0}^2}$ vanishes. 
\end{proof}

\begin{proposition}\label{prop:coeffdegree}
 For $\nu$ in the interior of $\mathcal{R}(\varphi\oplus\Psi)$ and $0\leq i\leq 2$ there exists a nonzero maximal minor of $\Theta_\nu$ of degree $N_i$ in the coefficients of $F_i$, where $N_i$ is given in Proposition~\ref{prop:degrees}.
\end{proposition}
\begin{proof}
Without loss of generality we assume $i=0$. Choose a specialization $F=(f_0,f_1,f_2)$ such that $F^{\rm sat}=G^{\rm sat}$ and such that the ideal  $F'=(f_1,f_2)$ has height two. In this case the variety $V(F':_RG)$ has degree  
$$\deg(F':_RG)=\deg(F')-\deg(G)= a_{1}b_{2} + b_{1}a_{2} - \sum_{i=1}^{p}e_{i}=N_0$$
 Denote by $\psi_{12}$ the submatrix of $\psi$ consisting of the columns corresponding to the coefficients of $f_1, f_2$. Since $F'\subseteq F':_R G$ we deduce $\het(F':_RG)\geq \het(F')=2$. In view of Lemma~\ref{lem:phipsi} and Lemma~\ref{lem:homology} we conclude that $F':_R G = I_2(\varphi \oplus \psi_{12})$ and $EN(\psi\oplus \psi_{12})$ is a resolution of $R/(F':_RG)$. Moreover by Corollary~\ref{cor:virtualR} since $F^{\rm sat}=G^{\rm sat}$  it follows that $EN(\varphi\oplus\psi)$ is a virtual projective resolution for $R/I_2(\varphi \oplus \psi)$.

Let  $ \mathcal{R}'=\cR(\varphi\oplus\psi_{12})$ denote the region specified by Proposition~\ref{prop:MSbound}, which is contained in the  weak regularity region of $R/I_n(\varphi\oplus \psi_{12})$ and let $\mathcal{R}=\cR(\varphi\oplus\psi)$ be the corresponding region for $R/I_n(\varphi\oplus \psi)$. Using  Remark~\ref{rem:subcomplex}, since  $\nu \in \mathcal{R}$  it follows that $\nu \in \mathcal{R}'$, hence $R/I_n(\varphi\oplus\psi_{12})$ is also $\nu$-regular. By \cite[Corollary 2.15.]{MaclaganSmith2} we deduce that $H_{R/(F':_RG)}(\mathbf{\nu})=N_0$. Therefore 
  \[
\dim_k \left(I_n\left(\varphi\oplus\psi_{12}\right)\right)_\nu =\dim_k (F':_R G)_\nu=\dim_k R_\nu-N_0.
\]
Denote by $\theta_{12}$ the matrix corresponding to the Eagon-Northcott complex of $\varphi\oplus\psi_{12}$.  Following the discussion before Proposition~\ref{pro:multiple}, the image of this matrix, $\left(I_n\left(\varphi\oplus\psi_{12}\right)\right)_\nu$, is the vector space 
\[
{\rm Span}_k\left\{m\cdot \Delta_{i_1,\ldots,i_q}: \mbox{ none of the columns $i_1,\ldots,i_q$ involve the column of coefficients of $F_0$ }\right\}.
\]
Hence we can choose exactly $\dim_k R_\nu-N_0$ columns in the matrix $\theta_\nu$ that are independent and do not involve the coefficients of $F_0$, therefore the same is true for the matrix $\Theta_\nu$. Denote the submatrix consisting of these columns by $\Theta_{\nu,F_1,F_2}$. Next, by Proposition~\ref{prop:surjective} it follows that the map $\theta_\nu$ is surjective and thus its image has dimension $\dim_k R_\nu$. Thus the vector space 
\[
{\rm Span}_k\left\{m\cdot \Delta_{i_1,\ldots,i_q}: \mbox{ one column in ${i_1,\ldots,i_q}$ is a coefficient column of $F_0$}\right \}
\] 
has dimension $N_0$. Therefore there exists $N_0$ linearly independent columns in $\theta_\nu$ that only involve the coefficients of $F_0$ and the same is true for $\Theta_\nu$. Denote the submatrix given by these columns matrix by $\Theta_{\nu,F_0}$. The columns of  $\Theta_{\nu,F_1,F_2}$ together with the columns of $\Theta_{\nu,F_0}$ span a vector space of dimension $\dim_k R_\nu$, hence the maximal minor corresponding to these columns is a maximal non vanishing minor of $\Theta_\nu$. Furthermore, since the entries of $\Theta_\nu$ are linear in the coefficients of $F_0$, the determinant of this minor has degree $N_0$ in the coefficients of $F_0$, as desired.
\end{proof}

\begin{proposition}
\label{prop:gcd}
The greatest common divisor of the maximal minors of the matrix $\Theta_\nu$ is exactly $\Res_{\mathcal{G},\{(a_i,b_i)\}_{i=0}^2}$.
\end{proposition}
\begin{proof}
 Let $d$ be the greatest common divisor of the maximal minors of $\Theta_\nu$. Proposition~\ref{pro:multiple} implies that $d$ is a multiple of $\Res_{G,\{(a_i,b_i)\}_{i=0}^2}$.  However, Proposition~\ref{prop:coeffdegree} states  that the degree of $d$ in the coefficients of $F_0$ is less than or equal to $N_0$ and on the other hand Proposition~\ref{prop:degrees} implies that $\Res_{G,\{(a_i,b_i)\}_{i=0}^2}$ has degree $N_0$ in the coefficients of $F_0$. Therefore the degree of $d$ in the coefficients of $F_0$ is equal to $N_0$. The same argument for $i=1,2$ allows to conclude that $d=\Res_{G,\{(a_i,b_i)\}_{i=0}^2}$ since they have the same degree with respect to all sets of coefficients.
\end{proof}

Proposition~\ref{prop:gcd} gives a practical method to compute the residual resultant. Note that Lemma~\ref{lem:hypothesisgeneric} yields that the Eagon-Northcott complex gives a virtual projective resolution in this context.

\begin{alg}[Computation of the residual resultant]\hfill
\hfill \label{alg:residualres}

\noindent {\bf Input:} $G$ a locally complete intersection ideal with syzygy matrix $\varphi$, $\Psi$ as in equation \eqref{eq:FandPsi}.
\label{alg:resultant}
\begin{enumerate}
\item Pick $\nu$ in the interior of the regularity region $\cR(\varphi\oplus\Psi)$.
\item Compute the matrix $\Theta_\nu$ as explained before Proposition \ref{pro:multiple}.
\item Compute a maximal minor $\delta_i$ of degree $N_i$ in the coefficients of $F_i$ for $0\leq i \leq 2$.
\item Return $\gcd\left(\det(\delta_0),\det(\delta_1),\det(\delta_2)\right)$.
\end{enumerate}
\end{alg}

Examples illustrating this algorithm can be found in section \ref{s:examplesres}.

\begin{remark} \label{rm:det}
The computations in steps $(3)$ and $(4)$ in the above algorithm are computationally expensive. 
However we can replace these two steps by the computation of the determinant of 
the bidegree $\nu$ strand of the complex $EN(\varphi\oplus\Psi)$.
Briefly, the determinant of
a complex is an alternating product of minors of the matrices of the differentials in the complex. 
Theorem 34 in \cite{GKZ}[Appendix A] establishes an equality between the
 the $\gcd$ of the maximal minors of the first differential of a complex and the determinant of a complex under
certain hypotheses. Such hypotheses are satisfied for the complex $EN(\varphi\oplus\Psi)$ 
and therefore we can use determinants of complexes in this setting. 
We refer the reader to 
Appendix A in 
\cite{GKZ} for a detailed construction of the determinant of a complex. Although computing the determinant of a 
complex can also be
computationally expensive, by comparison it is faster than computing the $\gcd$ of the maximal minors.
\end{remark}

\section{\large Implicitization of tensor product surfaces}\label{s:implicitization}

We now describe the specific setting of interest for our paper. First we establish the relation between the residual resultant and the implicit equation of $\Lambda$ and immediately after we give explicit steps for its computation. 
Setting the coordinate ring of $\P^3$ to be $S=k[X,Y,Z,W]$ our goal is to find the equation $H\in S$ defining the algebraic variety 
\begin{eqnarray*}
\Lambda=\overline{\img(\lambda)}=\left\{[x:y:z:w]\in \mathbb{P}_k^3 : p_0w-p_3x=p_1w-p_3y=p_2w-p_3z=0\right\}=V(H)
\end{eqnarray*}
where  $\lambda: \Q \dashrightarrow \P^3$ is a rational map as described in the introduction.

Let $P=\langle p_0,p_1,p_2,p_3 \rangle$
be the ideal of $R$ generated by the polynomials that define the parameterization $\lambda$ and set  $T=R\otimes_kS$. We assume that the $p_i$ have no
common factors and that $P$ is a height two ideal in $R$ that defines a local complete intersection set of points. Let $G=P^{\mathrm{sat}}$ denote
the $B$-saturated ideal that defines the set of points in $\Q$ and
set $\mathcal{G}$ to be the sheaf of ideals on $\Q$ associated to $G$. Since $P^{\mathrm{sat}}=G$, the sheaf $\mathcal{G}(a,b)$ is generated
by its global sections $p_0,p_1,p_2,p_3$ on $\Q\setminus V(G)$.
We denote by $\pi:\widetilde{Q}\to \Q$ the blow-up of $\Q$ along $\mathcal{G}$, and
by $\tilde{p_i}$ the global section $\pi^{*}(p_i)$ of the sheaf $\widetilde{\mathcal{G}}_{(a,b)}$
for $i=0,1,2,3$. Since $\widetilde{\mathcal{G}}_{(a,b)}$ is an invertible sheaf on $\widetilde{Q}$ and 
$\tilde{p}_0,\tilde{p}_1,\tilde{p}_2,\tilde{p}_3$ are global sections that generate it,
we deduce that there is a morphism
\[
\widetilde{\lambda}:\widetilde{Q}\to\pp^3
\]
such that $\widetilde{\lambda}^{*}\mathcal{O}(1)\cong \widetilde{\mathcal{G}}_{(a,b)}$ and $\widetilde{\lambda}^{*}(x)=
\tilde{p}_0$,
$\widetilde{\lambda}^{*}(y)=\tilde{p}_1$, $\tilde{\lambda}^{*}(z)=\tilde{p}_2$, $\widetilde{\lambda}^{*}(w)
=\tilde{p}_3$ (\cite[Ch.II.7]{Hartshorne}). As $\widetilde{Q}$ is projective and irreducible, we have
$\widetilde{\lambda}_{*}(\widetilde{Q})=\deg(\widetilde{Q}/\Lambda)\Lambda$ where $\Lambda$ is the rational surface in
$\pp^3$  and $\deg(\widetilde{Q}/\Lambda)$ is the degree of the surjective map
$\widetilde{\lambda}: \widetilde{Q} \to \Lambda$. 

Let $\beta$ be the following regular map
\begin{align*}
\beta:U= \Q  \setminus V(G)&\longrightarrow \pp^{3} \\
[s:t]\times[u,v]&\mapsto (p_0:p_1:p_2:p_3).
\end{align*}
\begin{proposition}
\label{prop:degree}
The degree of $\Lambda$ divides \[2ab -\sum_{i=1}^p e_i\]
where $e_i$ is defined before Remark~\ref{rem:incidence} and it is equal to this number when $\beta$ is birational.
\end{proposition}
\begin{proof}
We have $\deg(\widetilde{\lambda}_{*}(\widetilde{Q}))=\deg(\widetilde{Q}/\Lambda)\cdot \deg(\Lambda)$.
Next, we compute $\deg(\widetilde{\lambda}_{*}(\widetilde{Q}))$ by
\[\deg(\widetilde{\lambda}_{*}(\widetilde{Q}))=\int_{\widetilde{Q}}c_1(\widetilde{\lambda}^{*}\mathcal{O}(1))^2
= \int_{\widetilde{Q}} c_1(\widetilde{\mathcal{G}}_{(a,b)})^2= 2ab -\sum_{i=1}^p e_i.\]
The last equality above follows from the same computation as in the proof of Proposition~\ref{prop:degrees}.
Thus $\deg(\widetilde{\lambda}_{*}(\widetilde{Q}))= 2ab -\sum_{i=1}^p e_i$, which proves the first part of the statement.

Now we consider the following diagram, where $E$ denotes the exceptional divisor of the
blow-up $\pi$,
\[
\xymatrix{
  \widetilde{Q} \setminus E \ar[r]^{\widetilde{\lambda}\mid_{\widetilde{Q} \setminus E}} \ar[d]_{\pi} & \pp^3      \\
  U=\Q\setminus V(G) \ar[r]^{\hspace{1.2cm}\beta} &\pp^3.
   }
\]
Since by construction $\widetilde{\lambda}$ is unique and since the vertical map is an isomorphism outside
the exceptional divisor, we deduce that
$\widetilde{\lambda}\mid_{\widetilde{Q} \setminus E }=\beta \circ \pi$
and hence $\deg(\widetilde{Q}/\Lambda)=\deg(U/\beta(U))$ which is one if $\pi$ is birational.
\end{proof}
The next proposition establishes the relation between  residual resultanst in $\Q$ and the 
implicitization problem for tensor product surfaces with basepoints.

\begin{proposition}\label{prop:implicit}
Suppose that $(a,b)\geq (k_i,l_i)$ for all $i$,   
$(a,b)\geq (k_{j_1}+1,l_{j_1})$ for some $ j_1 $, and 
$(a,b)\geq (k_{j_2},l_{j_2}+1)$ for some $ j_2$.
 Then 
\begin{align} \label{eq:ressurf}
\mathrm{Res}_{\mathcal{G},(a,b)}(p_0-Xp_3,p_1-Yp_3,p_2-Zp_3)=H(X,Y,Z,1)^{\deg(U/\beta(U))} 
\end{align}
with $\deg(U/\beta(U))=1$ if $\beta$ is birational.
\end{proposition}
\begin{proof}
The residual resultant is defined as a general resultant over the blow-up of $\Q$ along $\mathcal{G}$.
Let $\xi$ denote a point in $\widetilde{Q}\setminus V(\widetilde{p}_3)$ and let $\widetilde{W}$ denote
the variety
\[\{\xi\times (x,y,z)\mid \tilde{p}_0(\xi)-x\tilde{p}_3(\xi)=\tilde{p}_1(\xi)-y\tilde{p}_3(\xi)=\tilde{p}_2
(\xi)-z
\tilde{p}_3(\xi)=0\}.\]
Note that considering only points in $\widetilde{Q}\setminus V(\widetilde{p}_3)$ for the incidence
variety is not a restriction. Indeed if $\xi$ is such that $\tilde{p}_3(\xi)=0$, then for some
$i\in\{0,1,2\}$ we must have $\tilde{p}_i(\xi)\neq 0$ because 
$\tilde{p}_0,\tilde{p}_1,\tilde{p}_2,\tilde{p}_3$ generate the sheaf $\tilde{\mathcal{G}}(a,b)$
on $\widetilde{Q}$. Thus $\xi$ cannot be a solution of the system $\tilde{p}_0(\xi)-x\tilde{p}_3(\xi)=\tilde{p}_1(\xi)-y\tilde{p}_3(\xi)=\tilde{p}_2
(\xi)-\tilde{p}_3z =0$.
Consider the following diagram
\[
\xymatrix{
  \widetilde{W} \ar[r]^{\pi_2} \ar[d]_{\pi_1} & \pp^3   \setminus V(W)   \\
  \widetilde{Q}\setminus V(\widetilde{p}_3)\ar[r]^{\hspace{-0.5cm}\pi}& \Q\setminus V(G) \ar[u]_{\beta}.
   }
\]
The  cycle in $\pp^3$ that represents the residual resultant is exactly ${\pi_{2}}_{*}(\widetilde{W})$,
i.e. ${\pi_{2}}_{*}(\widetilde{W})=\deg(\widetilde{W}/\pi_2(\widetilde{W}))\Lambda$ (in the generic case we
have $\deg(\widetilde{W}/\pi_2(\widetilde{W}))=1$). As the blow-up $\pi$ is an isomorphism outside the exceptional divisor,
the equation that defines ${\pi_{2}}_{*}(\widetilde{W})$ vanishes if and only if the point $(x,y,z,1)\in \pp^3$ is in
$\Lambda$. We deduce that
\[\mathrm{Res}_{\mathcal{G},(a,b)}(p_0-Xp_3,p_1-Yp_3,p_2-Zp_3)=H(X,Y,Z,1)^{\deg(\widetilde{W}/\pi_2(\widetilde{W}))}.\]
Now the map $\beta\mid_{\Q\setminus V(p_3)}$ makes the above diagram commute, and since $\pi$ is birational,
we deduce that $\deg(\widetilde{W}/\pi_2(\widetilde{W}))= \deg(U/\beta(U))$.

\end{proof}

\begin{remark}
It follows from Proposition~\ref{prop:degrees} that in the case where
$\deg(f_{0})=\deg(f_{1})=\deg(f_{2})=(a,b)$ then $\mathrm{Res}_{\mathcal{G},\{(a,b)\}}$
has degree $2ab-\sum_{i=1}^{p}e_{i}$ in the coefficients of $f_0,f_1$ and $f_2$. 
Looking at the degrees of the polynomials in equation (\ref{eq:ressurf})  from Proposition~\ref{prop:implicit}, we deduce that 
$2ab-\sum_{i=1}^{p}e_{i}=\deg(U/\beta(U)) \cdot \deg(H)$, this yields an alternate proof of the first assertion in Proposition~\ref{prop:degree}. \end{remark}

Proposition~\ref{prop:implicit} establishes that the residual resultant of $\F=(F_0,F_1,F_2)=(p_0-Xp_3,
p_1-Yp_3,p_2-Zp_3)$ with respect to $\mathcal{G}$ computes the implicit equation $H=0$.
To use the methods  presented in section \ref{s:resultantcomputation} to compute the implicit equation of a tensor product surface via residual
resultants we assume the given parameterization  has a special form. To set up a parametrization $\lambda$ we start with a locally complete intersection ideal $G =\langle g_1,\ldots,g_n\rangle\subseteq R$ of height two with
 a Hilbert-Burch resolution and four bihomogenous
polynomials $p_0,p_1,p_2$, $p_3 \in R_{(a,b)}$ related by
\begin{equation} \label{eq:pH}
\begin{bmatrix} p_0 & p_1 &p_2 & p_3 \end{bmatrix}= \begin{bmatrix} g_1 & \cdots & g_n\end{bmatrix} [h_{ji}],\quad h_{ji}
\in R_{(a-k_i,b-l_i)} .
\end{equation}
Second, we assume  $P^{\mathrm{sat}}=G^{\mathrm{sat}}$. The importance of this assumption is clarified in the following Lemma~\ref{lem:Psat}
and guarantees that we can use Eagon-Northcott complex of $I_n(\varphi\oplus\Psi)$ to find suitable degrees in the regularity region
of $T/I_n(\varphi\oplus\Psi)$.

\begin{lemma}
 \label{lem:Psat}
Suppose that $G$ is a locally complete intersection ideal and $P^{\rm sat}=G^{\rm sat}$. Then the ideal $\F=(p_0-Xp_3,p_1-Yp_3,p_2-Zp_3)$ has the property that $F_\fp=G_\fp$ for any ideal $\fp\in \Spec R \setminus\cB$ with $\het(\fp)=2$.
\end{lemma}

\begin{proof}
Let $\fp\in \Spec(R)\setminus \cB$ be an ideal with $\h (\fp)=2$.  Since $P^{\rm sat}=G^{\rm sat}$ it follows that $P_\fp=G_\fp$ and since $G_\fp$ is a complete intersection it is furthermore the case that $\dim_{k(\fp)}P_\fp/P^2_\fp=2$. Now $\F_\fp\subseteq P_\fp$ is generated by 3 elements which are pairwise independent in $P_\fp/P_\fp^2$, thus the equality $\F_\fp=G_\fp$ follows.
\end{proof}

Using the relation from Equation~(\ref{eq:pH}), we can write the polynomials $F_0=p_0-Xp_3, F_1=p_1-Yp_3, F_2=p_2-Zp_3$ from
Proposition~\ref{prop:implicit} as 
\begin{equation}
\label{eq:Psiimplicitization}
\begin{bmatrix} F_0 & F_1 & F_2 \end{bmatrix}=\begin{bmatrix} p_0-Xp_3 & p_1-Yp_3 & p_2-Zp_3\end{bmatrix} = \begin{bmatrix} g_1 & \cdots & g_n\end{bmatrix} \Psi.
\end{equation}
Based on Algorithm~\ref{alg:resultant} we derive a version that is tailored to the implicitization problem.

\begin{alg}[Implicitization algorithm]
\label{alg:implicitization}
\hfill

\noindent {\bf Input:} $G$ a locally complete intersection ideal, $P$ as in equation \eqref{eq:pH} such that $P^{\rm sat}=G^{\rm sat}$.
\begin{enumerate}
\item Set $\Psi=\begin{bmatrix} h_{i0}-Xh_{i3} & h_{i0}-Yh_{i3} & h_{i2}-Zh_{i3}\end{bmatrix}_{1\leq i\leq n}$ as in equation \eqref{eq:Psiimplicitization}.
\item Pick $\nu$ in the interior of the regularity region $\cR(\varphi\oplus\Psi)$.
\item Compute the matrix $\Theta_\nu$ as explained before Proposition \ref{pro:multiple}.
\item Compute a maximal minor $\delta_i$ of degree $N_i$ in the coefficients of $F_i$ for $0\leq i \leq 2$.
\item Return $\gcd\left(\det(\delta_0),\det(\delta_1),\det(\delta_2)\right)$.
\end{enumerate}
\end{alg}

Examples illustrating this algorithm can be found in section \ref{s:examplesimp}.

\begin{remark} \label{rm:resultantfromsubmaximalminors}
If the hypothesis $P^{\rm sat}=G^{\rm sat}$ is not satisfied, Algorithm \ref{alg:implicitization} no longer applies since the presentation map $\Theta$ for $R/I_n(\varphi\oplus \Psi)$ described in Proposition~\ref{prop:surjective}  is no longer surjective when restricted to any bidegree. However, given a bidegree $\nu\in \N^2$, if the dimension of the cokernel of $\Theta_\nu$ is $c$, then the proof of Proposition \ref{pro:multiple} shows that the resultant divides the greatest common divisor of the generators of the $c$-th Fitting ideal of $\Theta_\nu$, i.e. the minors of size $(\dim_k R_\nu-c)\times(\dim_k R_\nu-c)$ for $\Theta_\nu$. This is illustrated in Example \ref{ex:imp2}.
\end{remark}

\begin{remark}
 As highlighted in Remark~\ref{rm:det}, steps $(3)$ and $(4)$ can be replaced by the computation of the determinant 
 of a complex. Determinants of complexes are also used in
syzygy approach methods for implicitization of triangular and tensor product surfaces, see for instance \cite{chardin, 
residualBuse, Botbol}. 
More importantly, in the context of implicitization it is sufficient to compute $\Theta_{\nu}$. The matrix $\Theta_{\nu}$
is known as an implicit matrix representation of the surface. Matrix representations are a useful alternative to implicit
equations to represent a surface. 
A detailed account of their use in Geometric Modeling is outlined by Bus\'e \cite{BuseMatrixRep}.
\end{remark}


\section{\large Examples}
\label{s:examples}

 \subsection{Examples of computing residual resultants}
 \label{s:examplesres}
 \begin{example}[Residual resultant of one reduced point]
  \label{ex:onepoint}
 We compute the residual resultant $\Res_{G,(1,1)}$, where $G=\langle s,v\rangle$ is the defining ideal
 of the reduced point $[0:1]\times [1:0]$ in $\Q$. Consider the system
 \[
 \begin{array}{c} F_{0}=(uc_{00}+vc_{01})s+(sc_{02}+tc_{03})v\\
 F_{1}= (uc_{10}+vc_{11})s+(sc_{12}+tc_{13})v\\
 F_{2}= (uc_{20}+vc_{21})s+(sc_{22}+tc_{23})v\end{array},
 \]
  and let $T=C\otimes R$, where $C=k[c_{ij}]$ is the ring of generic coefficients. 
  The ideal $G$ is a complete intersection and the matrix $\varphi\oplus\Psi$ is
 \[\varphi\oplus\Psi= \left( \begin{array}{cccc} -v & uc_{00}+vc_{01} & uc_{10}+vc_{11} &uc_{20}+vc_{21} \\
 s & sc_{02}+tc_{03} & sc_{12}+tc_{13} & sc_{22}+tc_{23} 
 \end{array} \right). \]
 To calculate $\Res_{G,(1,1)}$, we  find a bidegree $\nu$ as in Remark~\ref{rem:regularity} and compute the matrix $\Theta_{\nu}$. Let $J$ denote the ideal $I_2(\varphi\oplus\Psi)$. From Proposition \ref{prop:ENreg}, since the numerical parameters for this example are $(a,b)=(c,d)=(e,f)=(1,1)$ and $(e_1,f_1)=(1,0), (e_2,f_2)=(0,1)$ we obtain the estimate

\begin{eqnarray*}
\mathcal{R} &=&  \left(3a+c-e-\min_{1\leq i\leq j\leq n}(e_i+e_j),3b+d-f-\min_{1\leq i\leq j\leq n}(f_i+f_j)\right)+\St_{-3}+\N^2 \\
&=& (3,3)+St_{-3}+\N^2 \\
&=&
\left((3,0)+\N^2\right)\cup \left((2,1)+\N^2\right) \cup \left((1,2)+\N^2\right)\cup  \left((0,3)+\N^2\right) \subseteq \reg(R/J).
\end{eqnarray*}

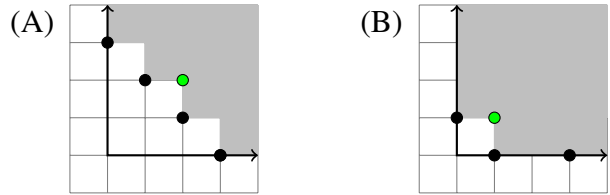
\begin{figure}[h!]
\label{fig:reg1}
\begin{tikzpicture}[scale=0.5]
\node (A)  at (-2,3.5){(A)};
   \draw[step=1cm,color=gray] (-1,-1) grid (4,4);
   \path [fill=lightgray] (0,3) rectangle (4,4);
   \path [fill=lightgray] (1,2) rectangle (4,3);
   \path [fill=lightgray] (2,1) rectangle (4,2);
   \path [fill=lightgray] (3,0) rectangle (4,1);
   \draw [<->] [thick] (0,4) -- (0,0) -- (4,0);
   \draw[fill] (0,3) circle [radius=0.15]; \draw[fill] (1,2) circle [radius=0.15];
   \draw[fill] (2,1) circle [radius=0.15]; \draw[fill] (3,0) circle [radius=0.15];
   \draw[fill=green] (2,2) circle [radius=0.15]; \draw[fill] (3,0) circle [radius=0.15];
\end{tikzpicture} \hspace{1cm}
\begin{tikzpicture}[scale=0.5]
\node (B)  at (-2,3.5){(B)};
   \draw[step=1cm,color=gray] (-1,-1) grid (4,4);
   \path [fill=lightgray] (0,1) rectangle (4,4);
   \path [fill=lightgray] (1,0) rectangle (4,1);
   \draw [<->] [thick] (0,4) -- (0,0) -- (4,0);
   \draw[fill] (0,1) circle [radius=0.15]; \draw[fill] (1,0) circle [radius=0.15];
   \draw[fill=green] (1,1) circle [radius=0.15]; \draw[fill] (3,0) circle [radius=0.15];
\end{tikzpicture}
\caption{
Example~\ref{ex:onepoint}, (A) regularity region $\mathcal{R}(\varphi\oplus\Psi)$, (B) strong regularity region.}
\end{figure}

We can choose any $\nu$ in the interior of the regularity region to set up $\Theta_{\nu}$. For $\nu=(2,2)$,
$\Theta_{(2,2)}$ is $9\times 24$ matrix. 
We can alternatively use the notion of strong regularity to find bidegrees such that $T/J$ is $\nu$-regular. 
  Computing the minimal free resolution for $J$ with {\em Macaulay2} \cite{M2} yields
  \[
   0 \rightarrow  T(-1,-2)\oplus T(-2,-1)\oplus T(-2,-2) \rightarrow  T(-1,-1)^2\oplus T(-1,-2)^3\oplus T(-2,-1)^3
  \]
  \[
\rightarrow  T(-1,-1)^6\rightarrow T \rightarrow T/J \rightarrow 0
  \]
  hence 
  \[
  \left((1,0)+\N^2\right) \cup \left((0,1)+\N^2\right) = \reg^s(T/J)\subseteq \reg(T/J).
  \]
   This means we can compute the determinant of the EN complex restricted to bidegree $(1,1)$ 
  to find the residual resultant of the system. 
  This yields the matrix $\Theta_{(1,1)}$ of size $6\times 12$.
 The residual resultant   is 
 \begin{align*}
\Res_{\mathcal{G},(1,1)}=-{c}_{03}{c}_{11}{c}_{20}-{c}_{03}{c}_{12}{c}_{20}+{c}_{01}{c}_{13}{c}_{20}+{c}_{02}{c}_{13}{c
      }_{20}+{c}_{03}{c}_{10}{c}_{21}-{c}_{00}{c}_{13}{c}_{21}+\\ {c}_{03}{c}_{10}{c}_{22}-{c}_{00}{c}_{1,
      3}{c}_{22}-{c}_{01}{c}_{10}{c}_{23}-{c}_{02}{c}_{10}{c}_{23}+{c}_{00}{c}_{11}{c}_{23}+{c}_{00}{a
      }_{12}{c}_{23}
 \end{align*}
For this example we can compute $\Res_{\mathcal{G},(1,1)}$ in a much simpler way. Indeed, we can 
 rewrite the system above as a
 linear system having three unkowns $su,sv,tv$. This system has the coefficient matrix 
 \[
 \rho=\left( \begin{array}{ccc} c_{00} & c_{01}+c_{02} &c_{03} \\ c_{10} & c_{11}+c_{12} & c_{13} \\ c_{20} & c_{21}+c_{22} & c_{23} \end{array} \right),\] 
hence  the system has a solution whenever this determinant is zero. Indeed, one can check that the displayed equation above gives $\Res_{G,(1,1)}=\det(\rho)$.
 \end{example}
 \begin{example}[Residual resultant of two complete intersection points]
 \label{ex:2}
 We compute the residual resultant $\Res_{\mathcal{G},(1,2)}$, where $G=\langle uv, s \rangle$ is a complete intersection defining a set of two reduced complete points in $\Q$ that lie on the same line in one of the rulings. Consider the system
        \[\begin{array}{ll}
        F_0 &= ( s {c}_{00}+t {c}_{01})uv+(u^{2} {c}_{02}+u v {c}_{03}+ v^{2} {c}_{04})s\\
      F_1 &= (s {c}_{10}+t {c}_{11}) uv+ (u^{2} {c}_{12}+u v {c}_{13}+v^{2} {c}_{14})s\\
      F_2 &= (s {c}_{20}+t {c}_{21})uv+   (u^{2} {c}_{22}+u v {c}_{23}+v^{2} {c}_{24})s.
      \end{array}\]
      According to Proposition~\ref{prop:degrees}, $\Res_{\mathcal{G},(1,2)}$
      is of degree $2$ in the coefficients of each $F_i$. We set up the  matrix
      \[
      \varphi\oplus \Psi=
      \bgroup\begin{pmatrix}{-s}&
      s {c}_{00}+t {c}_{01}&
      s {c}_{10}+t {c}_{11}&
      s {c}_{20}+t {c}_{21}\\
      u v&
      u^{2} {c}_{02}+u v {c}_{03}+v^{2} {c}_{04}&
      u^{2} {c}_{12}+u v {c}_{13}+v^{2} {c}_{14}&
      u^{2} {c}_{22}+u v {c}_{23}+v^{2} {c}_{24}
      \end{pmatrix}.\egroup\]
Let $J$ denote the ideal $I_2(\varphi\oplus\Psi)$. In a similar fashion as in Example~\ref{ex:onepoint}, we compute the regularity region for $R/J$ specified in 
Proposition~\ref{prop:ENreg} as illustrated in Figure~\ref{fig:reg2}. From this
region it follows that we may use $\nu=(1,6)$. The matrix $\Theta_{(1,6)}$ is of size $14\times 30$.
\begin{figure}[h!]

\begin{tikzpicture}[scale=0.5]
\node (A)  at (-2,6.5){(A)};
   \draw[step=1cm,color=gray] (-1,-1) grid (4,7);
   \path [fill=lightgray] (0,6) rectangle (4,7);
   \path [fill=lightgray] (1,5) rectangle (4,6);
   \path [fill=lightgray] (2,4) rectangle (4,5);
   \path [fill=lightgray] (3,3) rectangle (4,4);
   \draw [<->] [thick] (0,7) -- (0,0) -- (4,0);
   \draw[fill] (0,6) circle [radius=0.15]; \draw[fill] (1,5) circle [radius=0.15];
   \draw[fill] (2,4) circle [radius=0.15]; \draw[fill] (3,3) circle [radius=0.15];
   \draw[fill=green] (1,6) circle [radius=0.15]; 
\end{tikzpicture} \hspace{1cm}
\begin{tikzpicture}[scale=0.5]
\node (B)  at (-2,6.5){(B)};
   \draw[step=1cm,color=gray] (-1,-1) grid (4,7);
   \path [fill=lightgray] (0,3) rectangle (4,7);
   \path [fill=lightgray] (1,2) rectangle (4,3);
   \draw [<->] [thick] (0,7) -- (0,0) -- (4,0);
   \draw[fill] (0,3) circle [radius=0.15]; \draw[fill] (1,2) circle [radius=0.15];
   \draw[fill=green] (1,3) circle [radius=0.15]; 
\end{tikzpicture}
\caption{
Example~\ref{ex:2}, (A) Regularity region $\mathcal{R}(\varphi\oplus\Psi)$, (B) strong regularity region.} \label{fig:reg2}
\end{figure}
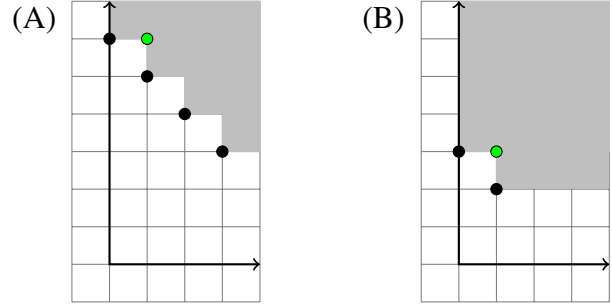 
\noindent The strong regularity region in this case is depicted in Figure~\ref{fig:reg2} (Right) and it is given by
$$\reg^s(T/J)=\left((1,2)+\Z_+^2\right) \cup \left((0,3)+\Z_+^2\right).$$
Estimating the regularity of $T/J$ using the strong regularity region allows the use of the bidegree $\nu=(1,3)$, for which the matrix $\Theta_{(1,3)}$ is an $8\times 12$ matrix. The polynomial
$\Res_{\mathcal{G},(1,2)}$ contains 141 terms.

         In this example, it is true that $\F:G=I_2(\varphi\oplus\Psi)$, but not for the reason given in the the hypothesis of Lemma~\ref{lem:phipsi} (2). The   Eagon-Northcott complex obtained from the matrix $\varphi\oplus\Psi$ is a virtual projective resolution for $I_2(\varphi\oplus\Psi)$, but it is not a resolution for this ideal since it is not exact. This is to be expected considering the proof of Lemma~\ref{lem:homology} because $\het (I_2(\varphi\oplus\Psi))=2$. 
However, this allows to estimate the regularity of $T/I_2(\varphi\oplus\Psi)$ using  Proposition \ref{prop:ENreg} as pictured in Figure~\ref{fig:reg2}.
   \end{example}

 \subsection{Examples of implicitization}
 \label{s:examplesimp}
 In this section we illustrate the techniques we developed  in the previous sections to compute the implicit
 equation of a map $\Q\to \pp^3$ defined by four bihomogeneous polynomials of bidegreee $(a,b)$.

 \begin{example}\label{ex:imp1}
 Let $I= \langle s,u \rangle \cap \langle t,v\rangle$ be  the ideal from Example~\ref{ex:twopoints}
 which defines two non-collinear points in $\Q$. This set is pictured below together with its Hilbert function.
\begin{center}
 \begin{picture}(70,50)

\put(0,10){\line(1,0){30}}
\put(0,20){\line(1,0){30}}

\put(7,35){$s$}
\put(17,35){$t$}

\put(30,18){$u$}
\put(30,8){$v$}

\put(10,20){\circle*{3}}
\put(20,10){\circle*{3}}

\put(10,0){\line(0,1){30}}
\put(20,0){\line(0,1){30}}
\put(45,20){\scalebox{0.6}{$H_{X}=\begin{array}{c||c|c|c|c} 
     & 0 & 1 & 2 & 3   \\ \hline \hline
    0& 1 & 2 & 2  & 2      \\ \hline
    1& 2  & 2  & 2  & 2    \\ \hline
    2& 2  & 2  & 2  & 2    \\ \hline
    3& 2 & 2 & 2 & 2 \\
   \end{array}$}}
\end{picture}
\end{center}
 Let
$G=\langle sv,tu\rangle$ and denote by $g_1,g_2$ the two generators of $G$. Here $G$ is a complete intersection with resolution
\[
\xymatrix{  0 \ar[r] 
&
R 
\ar[r]_{\scalebox{0.8}{ $\bgroup\begin{pmatrix}{-s v}\\
      t u\\
      \end{pmatrix}\egroup$
}} 
& 
R^2 \ar[r]_{\scalebox{0.8}{$\bgroup\begin{pmatrix} t u&  s v\\   \end{pmatrix}\egroup$}}
& 
G \ar[r]
&
0.
}
\]
Note that $G^{\mathrm{sat}}=I$, so, while $G$ is not saturated, however $V(G)=V(I)$ and therefore the complex displayed above 
is a Hilbert-Burch virtual resolution for $I$. Next we consider the  ideal 
$P=\langle p_0,p_1,p_2,p_3\rangle$ where 
$\begin{bmatrix} p_0 & p_1 & p_2 & p_3\end{bmatrix}=\begin{bmatrix} g_1 & g_2\end{bmatrix} h$ and $h$  
  is the $2\times 4$   matrix 
  \[h= \bgroup\begin{bmatrix}s&      t&      0&      0\\      0&     0&      s&      t\\
      \end{bmatrix}\egroup.\]
The bihomogeneous polynomials $p_0=s^2v,p_1=stv,p_2=stu,p_3=t^2u$ define a parameterization of a tensor product surface of bidegree $(2,1)$ with two basepoints given by $V(P)$. Note 
  that the homogeneous implicit equation for this surface is easily obtained and equal to
  $ YZ-XW=0$.  
Since the primary decomposition of the ideal $P$ is
  $P=(s^2,st,t^2)\cap (s,u)\cap (t,v)\cap (u,v)$, it follows that $P^\mathrm{sat}= G^{\mathrm{sat}}$. We obtain the matrix $\Psi$ 
  by writing 
  \[
\begin{bmatrix} p_0-Xp_3 & p_1-Yp_3 & p_2-Zp_3\end{bmatrix}
= \begin{bmatrix} g_1 & g_2\end{bmatrix}\underbrace{
\begin{bmatrix} h_{01} -h_{31}X & h_{11}- h_{31}Y & h_{21}- h_{31}Z \\ 
h_{02}-h_{32}X & h_{12}-h_{32}Y & h_{22}-h_{32}Z \end{bmatrix}}_\Psi.
\]
We note that the bidegree $(2,1)$ does not satisfy the inequality conditions in the hypotheses of 
Proposition~\ref{prop:tpres}. However, we can use the result in this proposition because we 
can find an open set such that
the sheaf $\widetilde{\mathcal{G}}(2,1)$ is very ample. To see this, set $U$ to be the open set described in the proof of 
Proposition~\ref{prop:tpres} with $g_{j_1}=g_1$ and $g_{j_2}=g_2$. It suffices to consider $x,y\in \Q$ and show that
the sections of $\mathcal{G}(2,1)$ separate points. Suppose that 
$I(x)=\langle l_1,h_2\rangle$,
$I(y)=\langle l_2,h_2\rangle$ with $l_i\in R_{(1,0)}$ and $h_i\in R_{(0,1)}$.   If $l_1$ is not a multiple of $l_2$ 
and $l_1=as+bt$, $a,b\in k$ then the form $svl_1\in P_{(2,1)}$  vanishes at $x$ and not at $y$. An analogous 
argument with $h_1,h_2$ shows that if $l_1$ is a multiple of $l_2$ we can find a form in $P_{(2,1)}$ that vanishes at
$x$ and not at $y$.   This shows that the pullbacks $\tilde{p}_0,\tilde{p}_1,\tilde{p}_2,\tilde{p}_3$ of $p_0,p_1,p_2,p_3$
to $\widetilde{Q}$ separate points. Following the proof in Proposition~\ref{prop:tpres}, we see that 
$\tilde{p}_0,\tilde{p}_1,\tilde{p}_2,\tilde{p}_3$ also separate tangents. Since $\widetilde{\mathcal{G}}(2,1)$ is generated by its global sections $\tilde{p}_0,\tilde{p}_1,\tilde{p}_2,\tilde{p}_3$ and $\widetilde{\mathcal{G}}(2,1)$ is very ample on an open subset, we conclude
the residual resultant $\mathrm{Res}_{\mathcal{G},(2,1)}$ exists and satisfies the same properties as in the conclusion 
of Proposition~\ref{prop:tpres}.

  To obtain the implicit equation using a residual 
  resultant we set up the matrix $\Theta_{\nu}$  for a bidegree $\nu$ according to Remark~\ref{rem:regularity}.
  On one hand we compute the regularity region of  $EN(\phi\oplus\Psi)$ 
  following Proposition~\ref{prop:ENreg}. On the other hand we compute the strong regularity region determined
  by a minimal free resolution of $T/J$. 
The regions found by these two methods and the shifts in the minimal free resolution of $R/J$ are displayed in Figure~\ref{fig:imp1}.
  
\begin{figure}[h!]
\label{fig:imp1}
\begin{tikzpicture}[scale=0.5]
\node (A)  at (-2,3.5){(A)};
   \draw[step=1cm,color=gray] (-1,-1) grid (5,4);
   \path [fill=lightgray] (1,1) rectangle (5,4);
   \path [fill=lightgray] (2,0) rectangle (5,1);
   \draw [<->] [thick] (0,4) -- (0,0) -- (5,0);
   \draw[fill] (1,1) circle [radius=0.15]; \draw[fill] (2,0) circle [radius=0.15];
\end{tikzpicture} \hspace{0.1cm}
\begin{tikzpicture}[scale=0.5]
\node (B)  at (-2,3.5){(B)};
   \draw[step=1cm,color=gray] (-1,-1) grid (5,4);
   \path [fill=lightgray] (1,1) rectangle (5,4);
   \path [fill=lightgray] (2,0) rectangle (5,1);
   \draw [<->] [thick] (0,4) -- (0,0) -- (5,0);
   \draw[fill] (1,1) circle [radius=0.15]; \draw[fill] (2,0) circle [radius=0.15];
   \draw[fill=green] (3,0) circle [radius=0.15];
   \node (C)  at (6,3.5){(C)};
\end{tikzpicture} 
\begin{tabular}[b]{|c|c|} \hline
$i$ & Shifts in homological degree $i$\\ \hline
1 & $(2,0),(2,1)$\\
2 & $(2,1),(3,0),(3,1)$\\
3 & $(3,1),(2,1),(3,0),(3,2)$\\
4 & $(3,2)$\\ \hline
\end{tabular} 

\caption{
Example~\ref{ex:imp1}, (A) regularity region from Proposition \ref{prop:ENreg}, (B) strong regularity, and (C) bigraded shifts of a minimal free resolution of $T/J$.} \label{fig:reg3}
\end{figure}
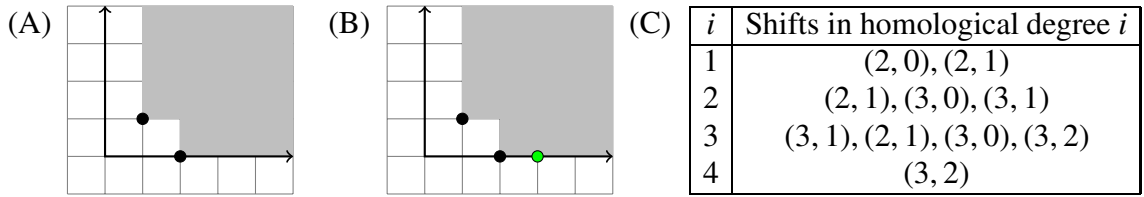 
Note that in this example, the two methods of estimating the regularity region for the module $T/I_n(\varphi\oplus\Psi)$ agree as shown in Figure~\ref{fig:imp1}.
 Now for $\nu=(3,0)$ one has 
\[ \Theta_{\nu} = \begin{bmatrix}
      0&0&1&0&0&0\\
      {-Y}&0&{-Z}&1&1&0\\
      X&{-Y}&0&{-Z}&{-Z}&1\\
      0&X&0&0&0&{-Z}\end{bmatrix}\]
   whence $I_4(\Theta_{\nu})=\langle YZ-X\rangle$ gives the implicit equation restricted to the affine set $W=1$. 
\end{example}

\begin{example} \label{ex:imp2}
Using the same setup as in Example~\ref{ex:imp1},  we change the entries of the matrix $h$ that determines
the parametrization ideal $P$. Set 
\[h= \bgroup\begin{bmatrix}s u&
       s v&
       0&
       t u+s v\\
       0&
       t u&
       s u&
       t v\\
       \end{bmatrix}\egroup, \]
       so $P= \langle s^{2}u\,v,t^{2}u^{2}+s^{2}v^{2},s\,t\,u^{2},s\,t\,u\,v+t^{2}u\,v+s^{2}v^{2}\rangle$.
       The generators of $P$ define a tensor product surface of bidegree $(2,2)$ with two basepoints.
       The support of $P$ and $G$ is the same, however the primary decomposition of $P$ reveals that the point corresponding to $(s,u)$ has multiplicity $2$ in the scheme defined by $P$. 

In this case we cannot use the Eagon-Northcott complex $EN(\varphi\oplus\Psi)$  to compute bidegrees
in the regularity region of $T/I_2(\varphi\oplus\Psi)$ because  the first homology module of  $EN(\varphi\oplus\Psi)$ is not $B$-torsion. In fact the first homology is a torsion module supported at the point with multiplicity 2 i.e.$ \langle s,u \rangle$. This shows the necessity of the hypothesis of Proposition~\ref{prop:virtualT}.

\begin{figure}[h!] 
\label{fig:imp2}
\begin{tikzpicture}[scale=0.5]
\node (A)  at (-2,4.5){(A)};
   \draw[step=1cm,color=gray] (-1,-1) grid (6,5);
   \path [fill=lightgray] (2,4) rectangle (6,5);
   \path [fill=lightgray] (3,3) rectangle (6,4);
   \path [fill=lightgray] (4,2) rectangle (6,3);
   \draw [<->] [thick] (0,5) -- (0,0) -- (6,0);
   \draw[fill] (2,4) circle [radius=0.15]; \draw[fill] (3,3) circle [radius=0.15];
   \draw[fill] (4,2) circle [radius=0.15]; 
   \draw[fill=green] (5,2) circle [radius=0.15]; 
   \node (B)  at (6.7,4.5){(B)};
\end{tikzpicture}
\hspace{-0.2cm}
\begin{tabular}[b]{|c|c|} \hline
$i$ & Shifts in homological degree $i$\\ \hline
1 & $(2,2)$\\
2 & $(4, 2), (2, 4), (4, 3), (3, 2), (3, 3), (2, 3)$\\
3 & $(4, 3), (4, 4), (4, 5),(3, 4), (3, 5), (5, 4), (2, 3), (4, 2), (5, 3)$\\
4 & $(3, 5), (5, 3), (3, 4), (5, 4), (5, 5)$\\ 
5 & $(4, 4), (4, 5), (5, 4), (5, 5)$\\ \hline
\end{tabular} 
\caption{ 
Example~\ref{ex:imp2}, (A) strong regularity and  (B) shifts in the resolution for $T/I_n(\varphi\oplus\Psi)$.} \label{fig:imp2}
\end{figure}
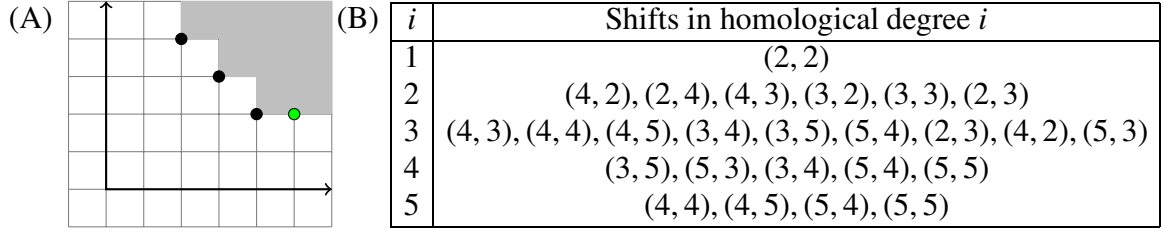
The free resolution of $T/I_2(\varphi\oplus\Psi)$  is 
$ 
0 \rightarrow  T^{5} \rightarrow T^{23}  \rightarrow T^{32}\rightarrow T^{19} \rightarrow T^{6} \rightarrow C
$ and the strong regularity region  
$
\reg^s\left(T/I_n(\varphi\oplus\Psi)\right)=\left((2,4)+\Z^2_+\right)\cup\left((3,3)+\Z^2_+\right)\cup \left((4,2)\Z^2_+\right)
$
is depicted in Figure~\ref{fig:imp2}. Therefore for $\nu=(5,2)$ the $18\times 24$  matrix $\Theta_{(5,2)}$ provides the implicit equation. 

Although $(3,2)$ is not the strong regularity region, we can use this bidegree to set up  a $12\times 12$ matrix $\Theta_{\nu}$ whose determinant vanishes, but that has an $11\times 11$ maximal minor whose determinant is a multiple of the implicit equation of the tensor product
surface. 
 \[\Theta_{(3,2),11\times 11}={\tiny\left({\begin{array}{ccccccccccc}
      0&0&0&0&0&0&0&0&1&0&0\\
      {-1}&0&0&0&0&0&0&0&{-X}&0&-Y+1\\
      X&0&Y-1&0&0&0&Z&0&0&0&0\\
      0&0&0&0&1&0&{-1}&0&{-X}&1&{-Y}\\
      X&{-1}&Y&0&-X-Y&0&Z&0&{-Z}&{-X}&Z\\
      0&X&0&Y-1&X&0&0&Z&0&0&{-Z}\\
      0&0&{-1}&0&{-X}&1&0&{-1}&0&{-X}&Z\\
      X&X&Y&Y&0&-X-Y&Z&Z&0&{-Z}&0\\
      0&0&0&0&0&X&0&0&0&0&0\\
      0&0&0&{-1}&0&{-X}&0&0&0&0&0\\
      0&X&0&Y&0&0&0&Z&0&0&0\\
      \end{array}}\right)}\]
 The implicit equation is the degree $5$ factor of \[X\cdot\left(X^{4}Y+X^{3}Y\,Z+X^{2}Y\,Z^{2}+X\,Y^{2}
 Z^{2}+X\,Y\,Z^{3}-X^{4}-2\,X^{2}Z^{2}-Z^{4}\right)\]
In this example the cokernel of $\Theta_{(3,2)}$ is 1-dimensional, and sum of the multiplicities of the basepoints is three, but there are two basepoints. This illustrates the observation made in Remark \ref{rm:resultantfromsubmaximalminors} that the residual resultant can be recovered as a divisor of the submaximal minors of $\Theta_{(3,2)}$ even if the base points in $P$ have higher multiplicity than the points in $G$.
 \end{example}
 
 \begin{example}  \label{ex:imp3}
 We continue with the setup from Example~\ref{ex:imp1} and change $h$ to
\[ h^\top =\begin{pmatrix}
       \frac{4}{9}\,s\,u+t\,u+s\,v+\frac{9}{5}\,t\,v&10\,s\,u+\frac{1}{2}\,t\,u+\frac{2}{3}\,s\,v+\frac{2}{3}\,t\,v\\
       \frac{1}{3}\,s\,u+\frac{10}{7}\,t\,u+\frac{9}{4}\,s\,v+\frac{2}{9}\,t\,v&\frac{8}{5}\,s\,u+\frac{1}{2}\,t\,u+\frac{5}{7}\,s\,v+\frac{2}{3}\,t\,v\\
       s\,u+\frac{4}{5}\,t\,u+s\,v+\frac{5}{8}\,t\,v&2\,s\,u+\frac{7}{3}\,t\,u+s\,v+\frac{9}{5}\,t\,v\\
       \frac{3}{5}\,s\,u+\frac{7}{3}\,t\,u+s\,v+8\,t\,v&\frac{4}{5}\,s\,u+\frac{7}{3}\,t\,u+\frac{3}{10}\,s\,v+\frac{7}{9}\,t\,v\end{pmatrix}.\]
 This choice of $h$ determines the ideal $P$ and a tensor product surface of bidegree $(2,2)$ with two basepoints 
 $V(P)$ and $P^{\mathrm{sat}}=G^{\mathrm{sat}}$. We use Proposition~\ref{prop:ENreg} to obtain the regularity region
 of $EN(\varphi\oplus\Psi)$
 depicted in Figure~\ref{fig:imp3}.
The resolution of  $T/I_2(\varphi\oplus\Psi) $ is  
\[
0 \rightarrow T^{119} \rightarrow T^{171} \rightarrow T^{71}\rightarrow T^{24} \rightarrow T^{6}\rightarrow T.
\]

\begin{figure}[h!]
\label{fig:imp3}
\begin{tikzpicture}[scale=0.5]
   \node (A)  at (-2,4.5){(A)};
   \draw[step=1cm,color=gray] (-1,-1) grid (5,5);
   \path [fill=lightgray] (1,4) rectangle (5,5);
   \path [fill=lightgray] (2,3) rectangle (5,4);
   \path [fill=lightgray] (3,2) rectangle (5,3);
   \path [fill=lightgray] (4,1) rectangle (5,2);
   \draw [<->] [thick] (0,5) -- (0,0) -- (5,0);
   \draw[fill] (1,4) circle [radius=0.15]; \draw[fill] (2,3) circle [radius=0.15];
   \draw[fill] (3,2) circle [radius=0.15]; \draw[fill] (4,1) circle [radius=0.15];
   \draw[fill=green] (3,3) circle [radius=0.15]; 
\end{tikzpicture} \hspace{1cm}
\begin{tikzpicture}[scale=0.5]
   \node (B)  at (-2,4.5){(B)};
   \draw[step=0.25cm,color=gray] (-1,-1) grid (5,5);
   \path [fill=lightgray] (3.5,3.25) rectangle (5,5);
   \draw [<->] [thick] (0,5) -- (0,0) -- (5,0);
   \draw[fill] (3.5,3.25) circle [radius=0.15]; 
\end{tikzpicture}
\caption{
Example~\ref{ex:imp3}, (A) regularity region of $EN(\varphi\oplus\Psi)$ and (B) strong regularity region with corner $(14,11)$.}
\end{figure}
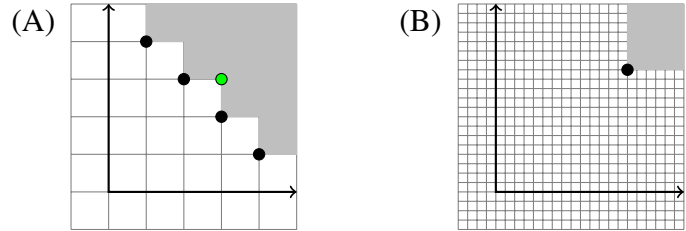 
The strong regularity region for this example is considerably worse than the regularity provided by Proposition \ref{prop:ENreg}.
For $\nu=(3,3)$, $\Theta_{\nu}$ is a matrix of size $16\times 24$. Although the point $\nu=(2,3)$ is not in the interior of the regularity
regions in Figure~\ref{fig:imp3}, the matrix $\Theta_\nu$ provides a $12\times 12$ determinental representation for the implicit 
equation of the surface.
 
 \end{example}



\section*{Acknowledgements}
Computations using the software system {\em Macaulay2} \cite{M2} were crucial for the development of this paper. We thank Laurent Bus\'e for suggesting that we explore residual resultants on $\Q$. The first author was supported by the
Deutsche Forschungsgemeinschaft(DFG, German Research Foundation )-314838170, GRK2297 MathCoRe.
The second author was supported by NSF grant DMS--1601024 and EpSCOR award OIA--1557417.

\bibliographystyle{amsalpha}	
\bibliography{refs}

\end{document}